\documentclass[hidelinks,onefignum,onetabnum]{siamart251216}

\usepackage[T1]{fontenc}
\usepackage[utf8]{inputenc}
\usepackage{amsfonts,amssymb,mathtools,bm}
\usepackage{graphicx,booktabs}
\usepackage{cite}
\usepackage{microtype}
\usepackage{enumitem,placeins,float}
\usepackage{algpseudocode}
\hypersetup{hidelinks,hypertexnames=false}
\allowdisplaybreaks

\newsiamthm{assumption}{Assumption}
\newsiamremark{remark}{Remark}

\newcommand{\R}{\mathbb{R}}
\newcommand{\Th}{\mathcal{T}_h}
\newcommand{\Ehzero}{\mathcal{E}_h^0}
\newcommand{\Z}{\bm Z}
\newcommand{\Zh}{\bm Z_h}
\newcommand{\Uh}{\bm U_h}
\newcommand{\vh}{\bm v_h}
\newcommand{\wh}{\bm w_h}
\newcommand{\zh}{\bm z_h}
\newcommand{\vbar}{\overline{\bm v}}
\newcommand{\Tbar}{\overline T}
\newcommand{\qbar}{\overline q}
\newcommand{\pbar}{\overline p}
\newcommand{\tnorm}[1]{\left\lVert #1\right\rVert_{V,h}}
\newcommand{\norm}[2][]{\left\lVert #2\right\rVert_{#1}}
\newcommand{\ip}[2]{\left(#1,#2\right)}
\newcommand{\dualnorm}[1]{\left\lVert #1\right\rVert_{\Zh'}}
\newcommand{\exampleheading}{%
  \refstepcounter{subsection}%
  \subsection*{Example~\thesubsection}%
  \addcontentsline{toc}{subsection}{Example~\thesubsection}%
}

\headers{Divergence-Free Convection-Cooling Control}{T. Wang and Z. Zhou}
\title{A Divergence-Free Finite Element Discretization of a Convection-Cooling Optimal Control Problem%
\thanks{This work was supported by the National Natural Science Foundation of China under Grant Nos.~12471385 and 11971276.}}

\author{Tongxin Wang\thanks{School of Mathematics and Statistics, Shandong Normal University, Jinan 250358, Shandong, China (Zhaojie Zhou, corresponding author: \email{zhouzhaojie@sdnu.edu.cn}).}
\and Zhaojie Zhou\footnotemark[2]}

\ifpdf
\hypersetup{
  pdftitle={A Divergence-Free Finite Element Discretization of a Convection-Cooling Optimal Control Problem},
  pdfauthor={Tongxin Wang and Zhaojie Zhou}
}
\fi

\begin{document}
\maketitle

\begin{abstract}
In this paper, we develop and analyze a finite element discretization for a bilinear optimal control problem governed by a stationary convection-diffusion equation, in which an incompressible velocity field acts as the control through the convective term. Continuous piecewise linear finite elements are employed to approximate the state and adjoint variables, while the Stokes-type optimality equation for the velocity control is discretized by a low-order, pointwise divergence-free ($P_1^c\oplus RT_0-P_0$) finite element. To accommodate the nonconforming velocity space, we construct a discrete cost functional by replacing the continuous $H^1$-seminorm regularization with a suitable discrete velocity bilinear form and derive the corresponding discrete first-order optimality system. Since the bilinear coupling between the velocity control and the state variable renders the underlying optimization problem nonconvex, we establish local strong convexity of the discrete reduced cost functional under a second-order sufficient optimality condition at an isolated strict local minimizer of the continuous problem. Furthermore, by exploiting a divergence-preserving interpolation operator, auxiliary state and adjoint error estimates, and a stability relation between the discrete reduced-gradient residual and the control error, we derive a priori error estimates for the state, adjoint state and velocity control. Numerical experiments are presented to demonstrate the effectiveness of the optimized incompressible velocity in redistributing the temperature field and to verify the predicted convergence rates.
\end{abstract}

\begin{keywords}
bilinear optimal control, convection-cooling, divergence-free finite element methods,
pressure robustness, a priori error estimates
\end{keywords}

\begin{MSCcodes}
49M25, 65N30, 65N15, 49K20
\end{MSCcodes}

\section{Introduction}

Optimal control problems governed by partial differential equations (PDEs) arise in a wide range of applications, including fluid dynamics, heat transfer, chemical transport, and environmental engineering. In many conventional PDE-constrained optimization problems, the control enters the governing equation as an additive source term or through boundary conditions \cite{GunzburgerHouSvobodny1991,RoschVexler2006,CasasMateosRaymond2007,Hinze2005,DedeQuarteroni2005,GongEtAl2018}. However, in certain applications, it is more natural to manipulate the transport mechanism itself rather than the external forcing. A representative example is the enhancement of heat transfer through active control of an incompressible velocity field. In this setting, the velocity acts as a multiplicative control through the convection term, leading to a bilinear optimal control problem \cite{HeHuMu2021-en,HuLiuWang2022}.

Bilinear, or multiplicative, optimal control problems arise when the control enters the governing PDE or its boundary conditions as a coefficient rather than as an additive source term.
The common classes are problems in which the control acts through a reaction coefficient \cite{KronerVexler2009,ShakyaSinha2020,AllendesCampanaOtarola2026,CasasWachsmuthWachsmuth2018,BersetcheEtAl2024}, a convection coefficient \cite{HeHuMu2021-en,GlowinskiEtAl2022} or Robin boundary coefficient \cite{Winkler2020,WangZhou2025}.
Although these models share a nonlinear control-to-state map and a generally nonconvex reduced objective, the optimality systems and numerical approximation techniques may differ substantially.

In this paper, we consider a stationary convection--diffusion bilinear optimal control problem of the form
\begin{equation*}
\begin{aligned}
  \underset{\substack{\bm v\in[H_0^1(\Omega)]^d\\
                       \nabla\cdot\bm v=0}}
             {\operatorname{\min}}\quad
  J(T,\bm v)
  &:=
  \frac12\norm[0]{T-\langle T\rangle}^2
  +\frac{\gamma}{2}\norm[0]{\nabla\bm v}^2,
  \\
  \text{subject to}\qquad
  -\kappa\Delta T+\bm v\cdot\nabla T
  &=f \quad\text{in }\Omega,
  \qquad
  T=0 \quad\text{on }\partial\Omega.
\end{aligned}
\end{equation*}
Here $\kappa>0$ is the thermal diffusivity, $\gamma>0$ is the
regularization parameter, $f$ is a prescribed heat source, and
$\langle T\rangle=|\Omega|^{-1}\int_\Omega T\,\mathrm dx$.
The objective is to optimize the distribution of the state while penalizing the magnitude of the velocity gradient. This formulation is particularly relevant to convection-cooling applications, where an appropriately designed incompressible flow can redistribute temperature and improve cooling performance. Unlike a standard distributed control, the velocity enters the differential operator directly and satisfies the incompressible condition. Consequently, the resulting optimization problem is generally nonconvex.

For bilinear optimal problems in which the control acts through a reaction coefficient or a Robin boundary coefficient extensive research work has already been carried out.
Yang, Chang, and Liu~\cite{YangChangLiu2008} obtained a priori error estimates and superconvergence results for a bilinear optimal control problem, using continuous piecewise linear approximations for the state and adjoint and piecewise constant approximations for the control. Kr\"oner and Vexler~\cite{KronerVexler2009} established a priori finite element error estimates for bilinear  elliptic optimal control problems, where three types of discretization for control variable were discussed. Borz\`i, Park, and Vallejos Lass~\cite{BorziParkVallejosLass2016} developed a multigrid algorithm for convection-diffusion problems with bilinear control.
  Fuica and Ot\'arola~\cite{FuicaOtarola2022} developed reliable and efficient a posteriori estimators for elliptic bilinear optimal control problems. Casas, Chrysafinos, and Mateos~\cite{CasasChrysafinosMateos2025} extended finite element error analysis to semilinear elliptic equations with reaction-coefficient controls, including convergence of discrete local solutions under second-order sufficient conditions. Wang, Wang, and Zhou~\cite{WangWangZhou2024} studied adaptive finite element approximation for bilinear optimal control governed by the fractional Laplacian. Winkler~\cite{Winkler2020} derived finite element error estimates for bilinear boundary control problems.

For bilinear optimal control problems with the control changing the transport operator through a term of the form $\boldsymbol v\cdot\nabla y$, there are still relatively few studies on numerical algorithms and numerical approximation. In the stationary convection-cooling setting, He, Hu, and Mu~\cite{HeHuMu2021-en} established the existence of optimal incompressible velocity controls, derived first- and second-order optimality conditions, and implemented the coupled optimality system using finite elements and Picard--Newton iterations. Glowinski, Song, Yuan, and Yue~\cite{GlowinskiEtAl2022} considered a time-dependent advection--reaction--diffusion system controlled by a divergence-free velocity field, establishing existence and first-order optimality results and designing a nested conjugate-gradient algorithm. Hu, Liu, and Wang~\cite{HuLiuWang2022} also investigated bilinear control of convection-cooling from open-loop and closed-loop perspectives. Unlike reaction-coefficient control, velocity control couples the control to a spatial derivative of the state and imposes an additional incompressibility constraint. The first-order optimality system consists of a coupled state equation, adjoint equation, and Stokes-type system for the velocity control.
The numerical treatment of the velocity control is closely related to finite element methods for incompressible flows.  In the context of optimal control of incompressible flows, Merdon and Wollner~\cite{MerdonWollner2023} analyzed the benefits of pressure-robust discretizations and demonstrated how restoring orthogonality to gradient forces can improve estimates for state and control variables. Their analysis concerns a different control setting, but highlights the potential significance of pressure-robust discretizations in PDE-constrained optimization.  Therefore, in the present work we aim to construct a divergence free finite element discretization and present rigorous error estimates for the discrete system. To the best of our knowledge, there is no error analysis for divergence-free finite element approximation of this kind of optimal control problem.

The main difficulties in the numerical treatment of the velocity-controlled
bilinear optimal control problem are threefold. First, the advection term
\(\bm v\cdot\nabla T\) is linear in the control \(\bm v\) and in the state
\(T\) separately, but not jointly; hence the control-to-state mapping is
nonlinear and the reduced cost functional is not convex. Second, the
control is subject to the divergence-free constraint \(\nabla\cdot\bm v=0\),
which is coupled with the state equation. In the first-order optimality
conditions this constraint turns the equation for the gradient into a
Stokes-type saddle-point problem, whose discretization must satisfy a discrete
inf--sup condition and whose velocity approximation should respect the
constraint as accurately as possible. Third, the velocity control enters the
state equation through \(\bm v\cdot\nabla T\). This derivative coupling amplifies the influence of the velocity
error on the state, adjoint, and control errors, and it interacts with the
pressure error whenever the velocity is obtained from a mixed discretization
of the Stokes system.

These difficulties motivate the use of an exactly divergence-free
discretization for the velocity control. Classical mixed methods, such as the
Taylor--Hood element used in \cite{HeHuMu2021-en}, only enforce the
divergence constraint weakly, so that \(\nabla\cdot\bm v_h\neq 0\) pointwise.
The resulting velocity error is then polluted by the pressure error, and the
discrete velocity may fail to be mass conservative, which is undesirable both
physically and numerically. In the context of the bilinear coupling
\(\bm v\cdot\nabla T\), such a pollution is particularly harmful: since the
control acts through the gradient of the state, an inaccurate or
pressure-dependent velocity directly contaminates the state and adjoint
approximations, and hence the reduced-gradient residual that controls the
error in the optimal control.

Motivated by these considerations, we develop and analyze a divergence-free finite element approximation for the stationary convection diffusion bilinear optimal control problem. Continuous piecewise linear finite elements are employed to approximate the state and adjoint variables, while the velocity-control equation is discretized by using a low-order, pointwise divergence-free $P_1^c\oplus RT_0-P_0$ element \cite{LiRui2021-en}. The main contributions of this work are summarized as follows.

\begin{enumerate}[label=\textup{(\arabic*)},leftmargin=2.6em,itemsep=0.5em]
\item \textbf{A divergence-free finite element formulation for the coupled optimality system.} Since the velocity space $P_1^c\oplus RT_0$ is $H(\mathrm{div})$-conforming but is not, in general, a subspace of $[H_0^1(\Omega)]^d$ owing to the possible tangential discontinuities of the Raviart--Thomas component, the continuous $H^1$-seminorm regularization cannot be directly restricted to the discrete velocity space. We construct a discrete cost functional by replacing the continuous $H^1$-seminorm regularization with a suitable discrete velocity bilinear form. The corresponding discrete first-order optimality system is derived, with the incompressibility constraint enforced pointwise by the chosen velocity--pressure finite element pair.

\item \textbf{Uniform local strong convexity of the discrete reduced cost functional.} Under a second-order sufficient optimality condition at an isolated strict local minimizer of the continuous problem, we establish local strong convexity of the discrete reduced cost functional, with a convexity constant independent of the mesh size. This property provides a fundamental stability mechanism for the subsequent analysis of the discrete optimality system.

\item \textbf{Optimal-order a priori error estimates for the coupled variables.} By combining a divergence-preserving interpolation operator, error estimates for auxiliary state and adjoint problems, and a stability relation between the discrete reduced-gradient residual and the velocity-control error, we establish optimal-order a priori error estimates for the state, adjoint, velocity control, and pressure. The analysis accounts for both the nonconformity of the velocity approximation and the nonlinear coupling arising from the bilinear control.

\item \textbf{Numerical verification and comparison with a classical mixed method.} Numerical experiments are performed to demonstrate the ability of the optimized incompressible velocity to redistribute the temperature field and to verify the predicted convergence rates. Comparisons with the classical MINI finite element method are also provided to assess the accuracy and numerical behavior of the proposed discretization.
\end{enumerate}

The remainder of this paper is organized as follows.
Section~2 reviews the continuous problem and its optimality conditions.
Section~3 develops the finite element discretization, its stability and
approximation properties, and the local error analysis.
Section~\ref{sec:numerics-en} presents the numerical experiments, followed
by concluding remarks.

\section{The continuous optimal control problem}

Let $\Omega\subset\R^d$, $d=2,3$, be a bounded connected Lipschitz domain with boundary
$\Gamma=\partial\Omega$. We define the following spaces:
\[
  L_0^2(\Omega)
  :=
  \left\{
  r\in L^2(\Omega):
  \int_\Omega r\,\mathrm dx=0,
  \right\}
\]
\[
  \Z
  :=
  \left\{
  \bm w\in[H_0^1(\Omega)]^d:
  \nabla\cdot\bm w=0
  \right\}
\]
and
\begin{align*}
  H(\operatorname{div};\Omega)
  &:=\left\{\bm w\in[L^2(\Omega)]^d:
           \nabla\cdot\bm w\in L^2(\Omega)\right\},\\
  H_0(\operatorname{div};\Omega)
  &:=\left\{\bm w\in H(\operatorname{div};\Omega):
           \bm w\cdot\bm n=0\text{ on }\Gamma\right\},
\end{align*}
where the divergence is understood in the distributional sense,
$\bm n$ is the outward unit normal to $\Gamma$, and the normal trace
vanishes in $H^{-1/2}(\Gamma)$.
For $\bm w\in[H_0^1(\Omega)]^d$ and $r\in L_0^2(\Omega)$, we write
$b(\bm w,r):=\ip{\nabla\cdot\bm w}{r}$.  The same definition applies to
$\bm w\in H_0(\operatorname{div};\Omega)$.
For $\phi\in L^2(\Omega)$, we set
\[
  \langle\phi\rangle
  :=
  \frac{1}{|\Omega|}
  \int_\Omega\phi\,\mathrm dx,
  \qquad
  D\phi:=\phi-\langle\phi\rangle.
\]
The operator $D$ is a self-adjoint orthogonal projection in $L^2(\Omega)$:
$D^*=D$ and $D^2=D$. Throughout the analysis, $\norm[0]{\cdot}$ denotes the $L^2(\Omega)$ norm.
For scalar or vector functions in $H_0^1(\Omega)$ we use the equivalent
energy norm
\(
  \norm[1]{w}:=\norm[0]{\nabla w}.
\)
All constants denoted by $C$ may change from line to line, but are independent
of the mesh size $h$ unless an additional dependence is stated explicitly.

For $\bm w\in[H_0^1(\Omega)]^d$ and $\eta,\phi\in H_0^1(\Omega)$, we define

\[
  c(\bm w;\eta,\phi)
  :=
  \int_\Omega(\bm w\cdot\nabla\eta)\phi\,\mathrm dx.
\]
For $\bm w\in\Z$, integration by parts establishes the skew-symmetry relations

\[
  c(\bm w;\eta,\phi)
  =
  -c(\bm w;\phi,\eta),
~
  c(\bm w;\eta,\eta)=0.
\]
For a control
$\bm v\in\Z$, the state $T=T(\bm v)$ solves
\begin{equation}
\label{eq:state-weak-en}
  \kappa\ip{\nabla T}{\nabla\phi}
  +
  c(\bm v;T,\phi)
  =
  \ip{f}{\phi}
  \qquad
  \forall\phi\in H_0^1(\Omega).
\end{equation}

The state and adjoint estimates and the derivative formulas below
follow the analysis in \cite{HeHuMu2021-en}.

\begin{lemma}
\label{lem:continuous-state-wellposed-en}
 For every $\bm v\in\Z$ and $f\in L^2(\Omega)$, the problem
\eqref{eq:state-weak-en} has a unique solution $T=T(\bm v)$, and
\begin{equation}
\label{eq:state-stability-en}
  \norm[1]{T}
  \leq
  C\kappa^{-1}\norm[0]{f}.
\end{equation}
\end{lemma}

The above well-posedness result defines the control-to-state map
\[
  S:\Z\longrightarrow H_0^1(\Omega),
  \qquad S(\bm v):=T(\bm v).
\]
Then the reduced optimal control problem is given below
\begin{equation}
\label{eq:optimization-en}
  \min_{\bm v\in\Z}j(\bm v):=J(S(\bm v),\bm v).
\end{equation}

\begin{proposition}
\label{prop:state-map-C2-en}
The map $S:\Z\to H_0^1(\Omega)$ is twice continuously Fr\'echet
differentiable.  For $\bm h\in\Z$, the derivative
$z=S'(\bm v)\bm h$ is the unique solution of
\begin{equation}
\label{eq:linearized-state-en}
  \kappa\ip{\nabla z}{\nabla\phi}
  +c(\bm v;z,\phi)
  =-c(\bm h;S(\bm v),\phi)
  \qquad\forall\phi\in H_0^1(\Omega),
\end{equation}
and
\begin{equation}
\label{eq:linearized-state-bound-en}
  \norm[1]{z}
  \leq C\kappa^{-1}\norm[1]{\bm h}\norm[1]{S(\bm v)}
  \leq C\kappa^{-2}\norm[0]{f}\norm[1]{\bm h}.
\end{equation}
For $\bm h_1,\bm h_2\in\Z$, we set
$z_i=S'(\bm v)\bm h_i$. The second derivative
$y=S''(\bm v)[\bm h_1,\bm h_2]$ is characterized by
\begin{equation}
\label{eq:second-linearized-state-en}
  \kappa\ip{\nabla y}{\nabla\phi}+c(\bm v;y,\phi)
  =-c(\bm h_1;z_2,\phi)-c(\bm h_2;z_1,\phi)
  \quad\forall\phi\in H_0^1(\Omega).
\end{equation}
\end{proposition}
\begin{proof}
The residual of \eqref{eq:state-weak-en} defines a smooth map from
$\Z\times H_0^1(\Omega)$ to $H^{-1}(\Omega)$, and its derivative with
respect to the state is an isomorphism by coercivity.
The implicit-function theorem and uniqueness therefore imply that $S$
is twice continuously Fr\'echet differentiable.
\end{proof}

To derive the first optimality system we introduce the adjoint state $q=q(\bm v)\in H_0^1(\Omega)$ defined by
\begin{equation}
\label{eq:adjoint-weak-en}
  \kappa\ip{\nabla q}{\nabla\psi}
  -
  c(\bm v;q,\psi)
  =
  \ip{DT}{\psi}
  \qquad
  \forall\psi\in H_0^1(\Omega).
\end{equation}

\begin{lemma}
\label{lem:continuous-adjoint-wellposed-en}
 For every $\bm v\in\Z$ and its corresponding state $T=T(\bm v)$,
problem \eqref{eq:adjoint-weak-en} admits a unique solution.  Moreover,
\begin{equation}
\label{eq:adjoint-stability-en}
  \norm[1]{q}
  \leq
  C\kappa^{-1}\norm[0]{DT}
  \leq
  C\kappa^{-2}\norm[0]{f}.
\end{equation}
\end{lemma}

Proposition~\ref{prop:state-map-C2-en} and the chain rule imply that
$j$ is twice continuously differentiable.  %
For a given $\bm v\in\Z$, we set $T=S(\bm v)$ and denote by $q=q(\bm v)$
the solution of \eqref{eq:adjoint-weak-en}.
 For every $\bm h\in\Z$, the reduced derivative is
\begin{equation}
\label{eq:continuous-gradient-en}
  j'(\bm v)\bm h
  =
  \gamma\ip{\nabla\bm v}{\nabla\bm h}
  -
  \ip{q}{\bm h \cdot\nabla T}.
\end{equation}

\begin{theorem}
\label{thm:continuous-first-order-en}
\cite[Theorem~3.2]{HeHuMu2021-en}.
Let $\vbar$ be a local minimizer of \eqref{eq:optimization-en}, and set
$\Tbar=S(\vbar)$ and $\qbar=q(\vbar)$.  Then there is a unique
$\pbar\in L_0^2(\Omega)$ such that
\begin{align*}
 \gamma\ip{\nabla\vbar}{\nabla\bm w}-b(\bm w,\pbar)
 &=\ip{\qbar\nabla\Tbar}{\bm w}
 &&\forall\bm w\in[H_0^1(\Omega)]^d,\\
 b(\vbar,\chi)&=0
 &&\forall\chi\in L_0^2(\Omega).
\end{align*}
Together with the state and adjoint equations, the following first optimality system holds
\begin{equation}
\label{eq:continuous-optimality-en}
\begin{cases}
  -\kappa\Delta\Tbar+\vbar\cdot\nabla\Tbar=f,\\
  -\kappa\Delta\qbar-\vbar\cdot\nabla\qbar
  =
  \Tbar-\langle\Tbar\rangle,\\
  -\gamma\Delta\vbar+\nabla\pbar=\qbar\nabla\Tbar,\\
  \nabla\cdot\vbar=0,
\end{cases}
\end{equation}
with $\Tbar=\qbar=0$ and $\vbar=\bm0$ on $\Gamma$.
\end{theorem}

The following second-order estimates and quadratic-growth argument
are based on \cite{HeHuMu2021-en}.
For $z=S'(\bm v)\bm h$, the second derivative of $j$ is
\begin{equation}
\label{eq:continuous-second-derivative-en}
  j''(\bm v)[\bm h,\bm h]
  =
  \norm[0]{Dz}^2
  +
  2\ip{z}{\bm h\cdot\nabla q}
  +
  \gamma\norm[0]{\nabla\bm h}^2.
\end{equation}
Because of the bilinear term $\bm v\cdot\nabla T$, the reduced functional is
generally nonconvex. We need the following second order condition.
\begin{assumption}[Second-order sufficient condition]
\label{ass:ssc-en}
Let $\vbar$ be a stationary point of $j$.  Assume that there exists
$\delta>0$ such that
\begin{equation}
\label{eq:ssc-en}
  j''(\vbar)[\bm h,\bm h]
  \geq
  \delta\norm[1]{\bm h}^2
  \qquad
  \forall\bm h\in\Z.
\end{equation}
\end{assumption}
\begin{lemma}
\label{lem:large-gamma-ssc-en}
 The condition
\(
  \gamma
  -
  C\kappa^{-4}\norm[0]{f}^2
  \geq
  \delta>0
\)
 ensures that Assumption~\ref{ass:ssc-en} holds with the same
lower-bound constant $\delta$.
\end{lemma}
In particular, Taylor's formula and the local Lipschitz continuity of
 the second derivative of $j$ imply
that there exist $\rho,c_0>0$ such that
\begin{equation}
\label{eq:continuous-quadratic-growth-en}
  j(\bm w)
  \geq
  j(\vbar)+c_0\norm[1]{\bm w-\vbar}^2
  \qquad
  \forall\bm w\in\Z,
  \quad
  \norm[1]{\bm w-\vbar}\leq\rho.
\end{equation}
Thus $\vbar$ is an isolated strict local minimizer.

\section{\texorpdfstring{Finite element discretization and error analysis}{Finite element discretization and error analysis}}

For the discretization, we use a conforming, shape-regular simplicial
mesh $\Th$ of $\Omega$ without hanging nodes. We denote the diameter of
$K\in\Th$ by $h_K$ and set $h=\max_{K\in\Th}h_K$.
We define
\[
  S_h
  :=
  \left\{
  \phi_h\in H_0^1(\Omega):
  \phi_h|_K\in P_1(K)
  \quad
  \forall K\in\Th
  \right\}.
\]
Both the   state and adjoint are approximated in $S_h$.
For the control variable, we use a velocity--pressure pair
for which the discrete incompressibility constraint ensures that
the discrete velocity is pointwise divergence-free.
Following \cite{LiRui2021-en}, we define
\[
  \bm V_h^1
  :=
  \left\{
  \bm w_h\in[H_0^1(\Omega)]^d:
  \bm w_h|_K\in[P_1(K)]^d\quad\forall K\in\Th
  \right\},
\]
\[
  RT_{0,h}^0
  :=
  \left\{
  \bm r_h\in H_0(\operatorname{div};\Omega):
  \bm r_h|_K\in[P_0(K)]^d+\bm xP_0(K)\quad\forall K\in\Th
  \right\},
\]
and
\[
  Q_h
  :=
  \left\{
  \chi_h\in L_0^2(\Omega):
  \chi_h|_K\in P_0(K)\quad\forall K\in\Th
  \right\}.
\]
The discrete velocity space is
\[
  \Uh:=\bm V_h^1\oplus RT_{0,h}^0.
\]
Every $\wh\in\Uh$ has a unique decomposition
\[
  \wh=\wh^1+\wh^R,
  \qquad
  \wh^1\in\bm V_h^1,
  \quad
  \wh^R\in RT_{0,h}^0.
\]
The direct-sum property follows from
$[H_0^1(\Omega)]^d\cap RT_{0,h}^0=\{\bm0\}$.

Using the divergence form $b$ defined in the continuous problem, we set
\begin{equation}
\label{eq:discrete-divfree-en}
  \Zh
  :=
  \left\{
  \wh\in\Uh:
  b(\wh,\chi_h)=0
  \quad
  \forall\chi_h\in Q_h
  \right\}.
\end{equation}

\begin{proposition}[Discrete mass conservation]
\label{prop:discrete-mass-conservation-en}
Every $\wh\in\Zh$ is pointwise divergence-free:
\[
  \nabla\cdot\wh\equiv0
  \qquad\text{in }\Omega.
\]
\end{proposition}

\begin{proof}
Since $\nabla\cdot\Uh\subset Q_h$, $\nabla\cdot\wh$ is piecewise constant for
$\wh\in\Uh$, and its global mean is zero because
$\wh\in H_0(\operatorname{div};\Omega)$.  Thus
$\nabla\cdot\wh\in Q_h$.  For $\wh\in\Zh$ one may choose
$\chi_h=\nabla\cdot\wh$ in \eqref{eq:discrete-divfree-en}, obtaining
\[
  \norm[0]{\nabla\cdot\wh}^2=0.
\]
Hence the divergence vanishes identically on every element, not merely in a
weak averaged sense.
\end{proof}

For
$\bm u_h=\bm u_h^1+\bm u_h^R$ and
$\wh=\wh^1+\wh^R$,  we define
\begin{equation}
\label{eq:ah-en}
  a_h(\bm u_h,\wh)
  :=
  \ip{\nabla\bm u_h^1}{\nabla\wh^1}
  +
  a_R(\bm u_h^R,\wh^R).
\end{equation}
Two admissible choices for the bilinear form $a_R$ are considered.
The first is the weighted mass form
\begin{equation}
\label{eq:aR0-en}
  a_R^0(\bm u_h^R,\wh^R)
  :=
  \sum_{K\in\Th}
  \alpha_K h_K^{-2}
  \ip{\bm u_h^R}{\wh^R}_K,
\end{equation}
where $0<\alpha_*\leq\alpha_K\leq\alpha^*$, with bounds independent of $h$.

The second choice is diagonal with respect to the standard $RT_0$ flux
basis. We denote by $\Ehzero$ the set of interior edges for $d=2$ or
interior faces for $d=3$, and we fix a unit normal $\bm n_e$ on each
$e\in\Ehzero$. The global basis functions
$\{\bm\Phi_e\}_{e\in\Ehzero}$ are characterized by
\[
  \int_{e'}\bm\Phi_e\cdot\bm n_{e'}\,\mathrm ds
  =\delta_{ee'}
  \qquad \forall e,e'\in\Ehzero.
\]
Writing
\[
  \bm u_h^R=\sum_{e\in\Ehzero}u_e\bm\Phi_e,
  \qquad
  \wh^R=\sum_{e\in\Ehzero}w_e\bm\Phi_e,
\]
we define the diagonal divergence form by
\begin{equation}
\label{eq:aRdiv-en}
  a_R^{\mathrm{div}}(\bm u_h^R,\wh^R)
  :=
  \sum_{K\in\Th}
  \sum_{\substack{e\in\Ehzero\\ e\subset\partial K}}
  \alpha_K u_e w_e
  \ip{\nabla\cdot\bm\Phi_e}
     {\nabla\cdot\bm\Phi_e}_K.
\end{equation}
Both choices satisfy the required bounds on shape-regular meshes
\cite[Lemma~3]{LiRui2021-en}:
\begin{equation}
\label{eq:aR-spectral-equivalence-en}
  c_R\sum_{K\in\Th}h_K^{-2}\norm[0,K]{\bm r_h}^2
  \leq
  a_R(\bm r_h,\bm r_h)
  \leq
  C_R\sum_{K\in\Th}h_K^{-2}\norm[0,K]{\bm r_h}^2
  \qquad
  \forall\bm r_h\in RT_{0,h}^0,
\end{equation}
where $c_R,C_R>0$ depend only on the shape-regularity and on the uniform
bounds for $\alpha_K$, but not on $h$.

We define the extended energy norm by
\begin{equation}
\label{eq:velocity-energy-norm-en}
  \tnorm{\bm w}^2
  :=
  a_h(\bm w,\bm w),
  \qquad
  \bm w\in[H_0^1(\Omega)]^d\oplus RT_{0,h}^0.
\end{equation}
The Li--Rui pair  \cite[Lemma~7]{LiRui2021-en}  satisfies the mesh-independent
inf--sup condition
\begin{equation*}
  \beta_0\norm[0]{\chi_h}
  \leq
  \sup_{\bm0\neq\wh\in\Uh}
  \frac{b(\wh,\chi_h)}{\tnorm{\wh}}
  \qquad
  \forall\chi_h\in Q_h,
\end{equation*}
where $\beta_0>0$ is independent of $h$.

We use the same notation $c$ for velocities in
$[H_0^1(\Omega)]^d\oplus RT_{0,h}^0$, with the integral evaluated elementwise.
If $\wh\in\Zh$ and $\eta_h,\phi_h\in S_h$, then normal continuity,
zero normal boundary trace, and elementwise divergence freedom imply
\begin{equation}
\label{eq:ch-skew-en}
  c(\wh;\eta_h,\phi_h)
  =
  -c(\wh;\phi_h,\eta_h),
  \qquad
  c(\wh;\eta_h,\eta_h)=0.
\end{equation}

\begin{lemma}
\label{lem:discrete-trilinear-en}
For $d=2,3$, the following estimates hold with a constant independent of
$h$:
\begin{align*}
  \norm[L^4(\Omega)]{\bm w}
  &\leq
  C\tnorm{\bm w},
  \\
  |c(\bm w;\eta,\phi)|
  &\leq
  C\tnorm{\bm w}\norm[1]{\eta}\norm[1]{\phi},
  \\
 \left|\ip{r}{\bm w\cdot\nabla\eta}\right|
  &\leq
  C\norm[1]{r}\norm[1]{\eta}\tnorm{\bm w}.
\end{align*}
Here $\bm w=\bm w^1+\bm w^R$ belongs to
$[H_0^1(\Omega)]^d\oplus RT_{0,h}^0$.
\end{lemma}

\begin{proof}

By Sobolev embedding, Poincar\'e inequality, and the definition of the
energy norm, we have
\[
  \norm[L^4(\Omega)]{\bm w^1}
  \leq
  C\norm[1]{\bm w^1}
  \leq
  C\tnorm{\bm w}.
\]
For the Raviart--Thomas part, the polynomial inverse inequality on a
shape-regular simplex implies
\[
  \norm[L^4(K)]{\bm w^R}
  \leq
  Ch_K^{-d/4}\norm[0,K]{\bm w^R}.
\]
Consequently, we obtain
\begin{align*}
  \norm[L^4(\Omega)]{\bm w^R}^4
  &=
  \sum_{K\in\Th}\norm[L^4(K)]{\bm w^R}^4\leq
  C\sum_{K\in\Th}h_K^{-d}\norm[0,K]{\bm w^R}^4
  \\
  &=
  C\sum_{K\in\Th}h_K^{4-d}
  \left(h_K^{-1}\norm[0,K]{\bm w^R}\right)^4.
\end{align*}
Since $d=2,3$ and $\Omega$ is fixed, $h_K^{4-d}\leq C_\Omega$.
Moreover, $\sum_K a_K^4\leq(\sum_Ka_K^2)^2$ for nonnegative $a_K$.
Taking $a_K=h_K^{-1}\norm[0,K]{\bm w^R}$ and using
\eqref{eq:aR-spectral-equivalence-en}, we obtain
\begin{align*}
  \norm[L^4(\Omega)]{\bm w^R}^4
  &\leq
  C\left(
    \sum_{K\in\Th}h_K^{-2}\norm[0,K]{\bm w^R}^2
  \right)^2\leq
  C a_R(\bm w^R,\bm w^R)^2
  \leq
  C\tnorm{\bm w}^4.
\end{align*}
The triangle inequality now proves
$\norm[L^4(\Omega)]{\bm w}\leq C\tnorm{\bm w}$.

For the first trilinear estimate, H\"older's inequality yields
\begin{align*}
  |c(\bm w;\eta,\phi)|
  &\leq
  \norm[L^4(\Omega)]{\bm w}
  \norm[0]{\nabla\eta}
  \norm[L^4(\Omega)]{\phi}
 \leq
  C\tnorm{\bm w}\norm[1]{\eta}\norm[1]{\phi}.
\end{align*}
Likewise, H\"older's inequality gives
\[
  \left|\ip{r}{\bm w\cdot\nabla\eta}\right|
  \leq
  \norm[L^4(\Omega)]{r}\norm[0]{\nabla\eta}
  \norm[L^4(\Omega)]{\bm w}
  \leq
  C\norm[1]{r}\norm[1]{\eta}\tnorm{\bm w}.
\]
\end{proof}

\subsection{Discrete optimality conditions}

The  discrete optimal control problem is defined as follows
\begin{equation}
\label{eq:full-discrete-optimization-en}
 \min_{(T_h,\vh)\in S_h\times\Zh}
 J_h(T_h,\vh)
 :=\frac12\norm[0]{DT_h}^2+\frac\gamma2a_h(\vh,\vh),
\end{equation}
subject to
\begin{equation}
\label{eq:discrete-state-en}
 \kappa\ip{\nabla T_h}{\nabla\phi_h}
 +c(\vh;T_h,\phi_h)=\ip f{\phi_h}
 \qquad\forall\phi_h\in S_h.
\end{equation}
Here the continuous $H^1$-seminorm
regularization is replaced by $a_h(\vh,\vh)$ defined in
\eqref{eq:ah-en}, since the discrete control space is not   nonconforming.
\begin{lemma}
\label{lem:discrete-state-wellposed-en}
For every $\vh\in\Zh$, equation \eqref{eq:discrete-state-en} has a unique
solution in $S_h$, and
\begin{equation}
\label{eq:discrete-state-stability-en}
 \norm[1]{T_h}\leq C\kappa^{-1}\norm[0]{f}.
\end{equation}
\end{lemma}

\begin{proof}
The bilinear form
\begin{equation}
\label{eq:discrete-convection-diffusion-form-en}
 A_h(\vh;\eta_h,\phi_h)
 :=\kappa\ip{\nabla\eta_h}{\nabla\phi_h}
   +c(\vh;\eta_h,\phi_h)
\end{equation}
is continuous by Lemma~\ref{lem:discrete-trilinear-en}.  Since
$\vh\in\Zh$, identity \eqref{eq:ch-skew-en} implies
$A_h(\vh;\phi_h,\phi_h)=\kappa\norm[1]{\phi_h}^2$.
The finite-dimensional Lax--Milgram theorem establishes existence and
uniqueness.  We test with $T_h$ and apply Poincar\'e's inequality to obtain
\[
 \kappa\norm[1]{T_h}^2=\ip f{T_h}
 \leq C\norm[0]{f}\norm[1]{T_h},\quad  \kappa\norm[1]{T_h} \leq C\norm[0]{f}.
\]
\end{proof}
As in the continuous case, we derive the following discrete first-order optimality system.
\begin{theorem}
\label{thm:discrete-first-order-en}
Let $(T_h,\vh)$ be a local solution of
\eqref{eq:full-discrete-optimization-en}.  Then there are unique
$q_h\in S_h$ and $p_h\in Q_h$ such that
\begin{subequations}
\label{eq:discrete-optimality-en}
\begin{align}
 \kappa\ip{\nabla T_h}{\nabla\phi_h}+c(\vh;T_h,\phi_h)
 &=\ip f{\phi_h},
 &&\forall\phi_h\in S_h,
 \label{eq:discrete-optimality-state-en}\\
 \kappa\ip{\nabla q_h}{\nabla\psi_h}-c(\vh;q_h,\psi_h)
 &=\ip{DT_h}{\psi_h},
 &&\forall\psi_h\in S_h,
 \label{eq:discrete-optimality-adjoint-en}\\
 \gamma a_h(\vh,\wh)-b(\wh,p_h)
 &=\ip{q_h\nabla T_h}{\wh},
 &&\forall\wh\in\Uh,
 \label{eq:discrete-optimality-control-en}\\
 b(\vh,\chi_h)&=0,
 &&\forall\chi_h\in Q_h.
 \label{eq:discrete-optimality-div-en}
\end{align}
\end{subequations}
\end{theorem}
\begin{proof}
The derivative of the state constraint with respect to $T_h$ is the
isomorphism in Lemma~\ref{lem:discrete-state-wellposed-en}.  The
finite-dimensional Lagrange multiplier theorem is therefore applicable.
We use the Lagrangian
\[
 \mathcal L_h(T_h,\vh,q_h)
 :=J_h(T_h,\vh)-A_h(\vh;T_h,q_h)+\ip f{q_h}.
\]
By taking variations in $q_h$ and $T_h$, we obtain the state and adjoint equations.
Variations $\zh\in\Zh$ in the control give
\[
 \gamma a_h(\vh,\zh)-\ip{q_h\nabla T_h}{\zh}=0
 \qquad\forall\zh\in\Zh.
\]
The functional on the left vanishes on the kernel $\Zh$ of the discrete
divergence.  The discrete inf--sup condition consequently supplies a unique
$p_h\in Q_h$, converting this identity into
\eqref{eq:discrete-optimality-control-en}.  Conversely, testing that equation
with $\zh\in\Zh$ recovers the reduced stationarity condition.
The last equation of the system also
implies $\nabla\cdot\vh=0$ pointwise by
Proposition~\ref{prop:discrete-mass-conservation-en}.
\end{proof}
\begin{lemma}
\label{lem:discrete-adjoint-wellposed-en}
For every $\vh\in\Zh$ and its associated state $T_h$, the adjoint equation
\eqref{eq:discrete-optimality-adjoint-en} has a unique solution $q_h\in S_h$,
which satisfies
\begin{equation}
\label{eq:discrete-adjoint-stability-en}
 \norm[1]{q_h}
 \leq C\kappa^{-1}\norm[0]{DT_h}
 \leq C\kappa^{-2}\norm[0]{f},
\end{equation}
where $C$ is independent of $h$, $\vh$, and $\kappa$.
\end{lemma}

\begin{proof}
The adjoint bilinear form is continuous by
Lemma~\ref{lem:discrete-trilinear-en} and coercive with constant $\kappa$
by \eqref{eq:ch-skew-en}.  Hence the adjoint equation is uniquely solvable.
We test with $\psi_h=q_h$ and apply Poincar\'e's inequality to obtain
\[
 \kappa\norm[1]{q_h}^2
 =\ip{DT_h}{q_h}
 \leq C\norm[0]{DT_h}\norm[1]{q_h}.
\]
The first bound follows by cancellation when $q_h\ne0$ and is immediate
otherwise.  The contraction $\norm[0]{DT_h}\leq\norm[0]{T_h}$,
Poincar\'e's inequality, and \eqref{eq:discrete-state-stability-en}
establish the second bound.
\end{proof}

The state equation defines the discrete control-to-state map
\[
 \mathcal S_h:\Zh\longrightarrow S_h,
 \qquad \mathcal S_h(\vh):=T_h.
\]
Its residual is polynomial in the finite element coefficients, and its
state derivative is invertible by
Lemma~\ref{lem:discrete-state-wellposed-en}.  The finite-dimensional
implicit-function theorem therefore shows that $\mathcal S_h$ is smooth.
The reduced functional associated with
\eqref{eq:full-discrete-optimization-en} is
\begin{equation}
\label{eq:discrete-functional-en}
 j_h(\vh):=J_h(\mathcal S_h(\vh),\vh)
 =\frac12\norm[0]{D\mathcal S_h(\vh)}^2
   +\frac\gamma2a_h(\vh,\vh),
 \qquad\vh\in\Zh.
\end{equation}

\begin{proposition}[Discrete reduced gradient]
\label{prop:discrete-gradient-en}
For every direction $\bm d_h\in\Zh$, the reduced derivative satisfies
\begin{equation}
\label{eq:discrete-gradient-en}
 j_h'(\vh)\bm d_h
 =\gamma a_h(\vh,\bm d_h)
  -\ip{q_h(\vh)}{\bm d_h\cdot\nabla\mathcal S_h(\vh)}.
\end{equation}
\end{proposition}

\begin{proof}
We set $z_h=\mathcal S_h'(\vh)\bm d_h$. By differentiating the discrete
state equation, we find that $z_h$ is uniquely determined by
\begin{equation}
\label{eq:discrete-linearized-state-en}
 A_h(\vh;z_h,\phi_h)
 =-c(\bm d_h;T_h,\phi_h)
 \qquad\forall\phi_h\in S_h.
\end{equation}
Using the chain rule, the adjoint equation, \eqref{eq:discrete-linearized-state-en},
and skew-symmetry, we obtain
\begin{align*}
 j_h'(\vh)\bm d_h
 &=\gamma a_h(\vh,\bm d_h)
   +\ip{T_h-\langle T_h\rangle}{z_h-\langle z_h\rangle}\\
 &=\gamma a_h(\vh,\bm d_h)
   +\ip{T_h-\langle T_h\rangle}{z_h}\\
 &=\gamma a_h(\vh,\bm d_h)+A_h(\vh;z_h,q_h)\\
 &=\gamma a_h(\vh,\bm d_h)-\ip{q_h\nabla T_h}{\bm d_h}.
\end{align*}
This proves \eqref{eq:discrete-gradient-en}.
Thus evaluation of the gradient requires
the state and adjoint, but not an explicit computation of $z_h$.
\end{proof}

For $\bm d_h\in\Zh$, we denote by
$z_h=\mathcal S_h'(\vh)\bm d_h$ the solution of
\eqref{eq:discrete-linearized-state-en}. We differentiate the state
equation twice and use the adjoint equation to eliminate the second state
derivative. We thereby obtain
\begin{equation}
\label{eq:discrete-second-derivative-en}
 j_h''(\vh)[\bm d_h,\bm d_h]
 =\norm[0]{Dz_h}^2
  +2\ip{z_h}{\bm d_h\cdot\nabla q_h}
  +\gamma\tnorm{\bm d_h}^2.
\end{equation}

\begin{proposition}
\label{prop:discrete-second-derivative-lower-bound-en}
For every $\vh,\bm d_h\in\Zh$, the following lower bound holds:
\begin{equation}
\label{eq:discrete-second-derivative-bound-en}
 j_h''(\vh)[\bm d_h,\bm d_h]
 \geq\left(\gamma-C\kappa^{-4}\norm[0]{f}^2\right)
 \tnorm{\bm d_h}^2.
\end{equation}
\end{proposition}

\begin{proof}
We test \eqref{eq:discrete-linearized-state-en} with $z_h$ and use
Lemma~\ref{lem:discrete-trilinear-en} and
\eqref{eq:discrete-state-stability-en} to obtain
\begin{equation}
\label{eq:discrete-linearized-state-bound-en}
 \norm[1]{z_h}\leq C\kappa^{-1}\tnorm{\bm d_h}\norm[1]{T_h}
 \leq C\kappa^{-2}\norm[0]{f}\tnorm{\bm d_h}.
\end{equation}
Together with the adjoint bound in \eqref{eq:discrete-adjoint-stability-en}, this implies
\[
 2\left|\ip{z_h}{\bm d_h\cdot\nabla q_h}\right|
 =2|c(\bm d_h;q_h,z_h)|
 \leq C\kappa^{-4}\norm[0]{f}^2\tnorm{\bm d_h}^2.
\]
Combining this estimate with \eqref{eq:discrete-second-derivative-en} proves
\eqref{eq:discrete-second-derivative-bound-en}.
\end{proof}
\begin{corollary}
\label{cor:discrete-large-gamma-en}
For sufficiently large $\gamma$,
Proposition~\ref{prop:discrete-second-derivative-lower-bound-en} implies
\[
 j_h''(\vh)[\bm d_h,\bm d_h]\geq\mu\tnorm{\bm d_h}^2
 \qquad\forall\vh,\bm d_h\in\Zh,
\]
where $\mu>0$ is independent of $h$.
\end{corollary}

\subsection{\texorpdfstring{ Some auxiliary results}{ Some auxiliary results}}

\begin{lemma}[Regularity of the optimality system]
\label{ass:regularity-en}
Let $\Omega$ be a bounded convex polygon or polyhedron, and let $f\in L^2(\Omega)$.
Every weak solution $(\Tbar,\qbar,\vbar,\pbar)$ of
\eqref{eq:continuous-optimality-en} satisfies
\[
 \Tbar,\qbar\in H^2(\Omega)\cap H_0^1(\Omega),\qquad
 \vbar\in[H^2(\Omega)\cap H_0^1(\Omega)]^d,
 \quad \pbar\in H^1(\Omega)\cap L_0^2(\Omega).
\]
\end{lemma}
The regularity of $\Tbar$ and $\qbar$ follows from the
convection--diffusion regularity results in
\cite{CasasMateosRosch2020}.
Since $H^2(\Omega)\hookrightarrow L^\infty(\Omega)$ for $d\leq3$,
the control equation has right-hand side
$\qbar\nabla\Tbar\in[L^2(\Omega)]^d$.
The regularity of $(\vbar,\pbar)$ then follows from the
Stokes regularity results on convex polygons and polyhedra
\cite{KelloggOsborn1976,Dauge1989}.

To derive the error estimates, we first introduce an interpolation operator
that preserves the divergence constraint.  Let $I_h$ be the continuous nodal
$P_1$ interpolant,
$\Pi_h^R$ the standard $RT_0$ interpolant, and $P_h$ the $L^2$ projection
onto $Q_h$. For
$\bm v\in[H^2(\Omega)\cap H_0^1(\Omega)]^d$,  we define
\begin{equation}
\label{eq:Pi-en}
  \Pi_h\bm v
  :=
  I_h\bm v
  +
  \Pi_h^R(\bm v-I_h\bm v).
\end{equation}
\begin{lemma}
\label{lem:divergence-preserving-interpolation-en}
The operator $\Pi_h$ defined by \eqref{eq:Pi-en} satisfies
\begin{equation}
\label{eq:commuting-en}
  \nabla\cdot\Pi_h\bm v
  =
  P_h(\nabla\cdot\bm v).
\end{equation}
Hence $\bm v\in\Z$ implies $\Pi_h\bm v\in\Zh$. If
$\bm v\in[H^2(\Omega)\cap H_0^1(\Omega)]^d$, then
\begin{equation}
\label{eq:velocity-interpolation-en}
  \norm[0]{\bm v-\Pi_h\bm v}
  +
  h\tnorm{\bm v-\Pi_h\bm v}
  \leq
  Ch^2|\bm v|_2.
\end{equation}
For $\phi\in H^2(\Omega)\cap H_0^1(\Omega)$ and
$p\in H^1(\Omega)\cap L_0^2(\Omega)$, the scalar and pressure approximations satisfy
\begin{align}
  \norm[0]{\phi-I_h\phi}
  +
  h\norm[1]{\phi-I_h\phi}
  &\leq
  Ch^2|\phi|_2,
  \label{eq:scalar-interpolation-en}\\
  \norm[0]{p-P_hp}
  &\leq
  Ch|p|_1.
  \label{eq:pressure-projection-en}
\end{align}
\end{lemma}
The commuting property and velocity estimate follow from
\cite[Lemma~4]{LiRui2021-en}. The spectral equivalence in
\eqref{eq:aR-spectral-equivalence-en} extends the velocity estimate to
the admissible forms $a_R$.
The scalar interpolation estimate follows from
\cite[Proposition~2.2.1]{BoffiBrezziFortin2013}.
The pressure projection estimate follows from the
elementwise Poincar\'e inequality.

Because $\Uh\not\subset[H_0^1(\Omega)]^d$, we define
\begin{equation}
\label{eq:consistency-en}
  \delta_h(\bm v,\wh)
  :=
  \ip{-\Delta\bm v}{\wh}
  -
  a_h(\bm v,\wh).
\end{equation}

\begin{lemma} \cite{LiRui2021-en}
\label{lem:nonconforming-consistency-en}
For every
$\bm v\in[H^2(\Omega)\cap H_0^1(\Omega)]^d$ and
$\wh=\wh^1+\wh^R\in\Uh$, the consistency residual satisfies
\[
  \delta_h(\bm v,\wh)
  =
  \ip{-\Delta\bm v}{\wh^R},
\]
and therefore
\begin{equation}
\label{eq:consistency-bound-en}
  |\delta_h(\bm v,\wh)|
  \leq
  Ch|\bm v|_2\tnorm{\wh}.
\end{equation}
\end{lemma}

For a continuous optimal solution with the regularity in
Lemma~\ref{ass:regularity-en}, the $H(\operatorname{div})$ Green formula
implies $\ip{\nabla\pbar}{\wh}=-b(\wh,\pbar)$ for every $\wh\in\Uh$;
the boundary term vanishes because $\wh\cdot\bm n=0$. We test the control
equation in \eqref{eq:continuous-optimality-en} and use
\eqref{eq:consistency-en} to obtain the identity
\begin{equation}
\label{eq:continuous-control-tested-en}
 \gamma a_h(\vbar,\wh)-b(\wh,\pbar)
 =\ip{\qbar\nabla\Tbar}{\wh}-\gamma\delta_h(\vbar,\wh)
 \qquad\forall\wh\in\Uh.
\end{equation}

We next separate approximation errors from the perturbations induced by the
discrete control.
For the continuous optimal solution
$(\Tbar,\qbar,\vbar,\pbar)$, we introduce
\[
 (\widehat T_h,\widehat q_h,\widehat{\bm v}_h,\widehat p_h)
 \in S_h\times S_h\times\Uh\times Q_h
\]
through the following three discrete problems.  The auxiliary state satisfies
\begin{equation}
\label{eq:aux-state-en}
  \kappa\ip{\nabla\widehat T_h}{\nabla\phi_h}
  +
  c(\vbar;\widehat T_h,\phi_h)
  =
  \ip{f}{\phi_h}
  \qquad
  \forall\phi_h\in S_h.
\end{equation}
The auxiliary adjoint is determined by
\begin{equation}
\label{eq:aux-adjoint-en}
  \kappa\ip{\nabla\widehat q_h}{\nabla\psi_h}
  -
  c(\vbar;\widehat q_h,\psi_h)
  =
  \ip{\Tbar-\langle\Tbar\rangle}{\psi_h}
  \qquad
  \forall\psi_h\in S_h.
\end{equation}
The auxiliary velocity--pressure pair satisfies
\begin{subequations}
\label{eq:aux-stokes-en}
\begin{align}
  \gamma a_h(\widehat{\bm v}_h,\wh)
  -
  b(\wh,\widehat p_h)
  &=
  \ip{\qbar\nabla\Tbar}{\wh},
  &&\forall\wh\in\Uh,\\
  b(\widehat{\bm v}_h,\chi_h)&=0,
  &&\forall\chi_h\in Q_h.
\end{align}
\end{subequations}
We set
\[
  M_{\vbar}
  :=
  1+\kappa^{-1}\norm[1]{\vbar}.
\]
\begin{lemma}
\label{lem:auxiliary-errors-en}
Under the regularity setting of Lemma~\ref{ass:regularity-en}, the auxiliary variables satisfy
\begin{align}
 & \norm[1]{\Tbar-\widehat T_h}
\leq
  ChM_{\vbar}|\Tbar|_2,
  \label{eq:aux-state-error-en}\\
  &\norm[1]{\qbar-\widehat q_h}
\leq
  ChM_{\vbar}|\qbar|_2,
  \label{eq:aux-adjoint-error-en}\\
  &\tnorm{\vbar-\widehat{\bm v}_h}
\leq
  Ch|\vbar|_2,
  \label{eq:aux-velocity-error-en}\\
  &\norm[0]{\pbar-\widehat p_h}
\leq
  Ch\bigl(|\pbar|_1+\gamma|\vbar|_2\bigr).
  \label{eq:aux-pressure-error-en}
\end{align}
\end{lemma}

\begin{proof}

We introduce the bilinear form
\[
  A_{\vbar}(\eta,\phi)
  :=
  \kappa\ip{\nabla\eta}{\nabla\phi}
  +
  c(\vbar;\eta,\phi).
\]
Galerkin orthogonality reads
$A_{\vbar}(\Tbar-\widehat T_h,\phi_h)=0$.
We decompose the error as
\[
  \Tbar-\widehat T_h
  =
  \eta_T+\theta_T,
  \quad
  \eta_T:=\Tbar-I_h\Tbar,
  \quad
  \theta_T:=I_h\Tbar-\widehat T_h.
\]
Since $c(\vbar;\theta_T,\theta_T)=0$, we have
\begin{align*}
  \kappa\norm[1]{\theta_T}^2
  &=
  A_{\vbar}(\theta_T,\theta_T)
  =
  -A_{\vbar}(\eta_T,\theta_T)
 \leq
  \left(
  \kappa+C\norm[1]{\vbar}
  \right)
  \norm[1]{\eta_T}\norm[1]{\theta_T}.
\end{align*}
We then apply the triangle inequality and the
interpolation estimate to obtain
\[
  \norm[1]{\Tbar-\widehat T_h}
  \leq
  CM_{\vbar}\norm[1]{\Tbar-I_h\Tbar}
  \leq
  ChM_{\vbar}|\Tbar|_2.
\]
The same argument applied to the adjoint form
\[
  A_{\vbar}^*(\eta,\psi)
  :=
  \kappa\ip{\nabla\eta}{\nabla\psi}
  -
  c(\vbar;\eta,\psi)
\]
establishes \eqref{eq:aux-adjoint-error-en}.

For the Stokes variables, the commuting property shows that
$\Pi_h\vbar\in\Zh$. The second equation in
\eqref{eq:aux-stokes-en} and the pointwise-divergence argument show that
$\widehat{\bm v}_h\in\Zh$. We set
\[
  \bm\xi_h
  :=
  \widehat{\bm v}_h-\Pi_h\vbar
  \in\Zh.
\]
We subtract \eqref{eq:continuous-control-tested-en} from
\eqref{eq:aux-stokes-en} and choose
$\wh=\bm\xi_h$. This gives
\[
  \gamma a_h(\widehat{\bm v}_h-\vbar,\bm\xi_h)
  -b(\bm\xi_h,\widehat p_h-\pbar)
  =\gamma\delta_h(\vbar,\bm\xi_h).
\]
The pressure term vanishes because $\bm\xi_h\in\Zh$ is pointwise
divergence-free by Proposition~\ref{prop:discrete-mass-conservation-en}.
After inserting
$\widehat{\bm v}_h-\vbar
=\bm\xi_h+(\Pi_h\vbar-\vbar)$ and dividing by $\gamma>0$, we obtain
\begin{align*}
  \tnorm{\bm\xi_h}^2
  &=a_h(\vbar-\Pi_h\vbar,\bm\xi_h)
    +\delta_h(\vbar,\bm\xi_h)
 \leq
  \left(
    \tnorm{\vbar-\Pi_h\vbar}+Ch|\vbar|_2
  \right)\tnorm{\bm\xi_h}.
\end{align*}
Equations \eqref{eq:velocity-interpolation-en} and
\eqref{eq:consistency-bound-en} show that
$\tnorm{\bm\xi_h}\leq Ch|\vbar|_2$; a triangle inequality proves
\eqref{eq:aux-velocity-error-en}.

Finally, we set
$\eta_p:=\widehat p_h-P_h\pbar\in Q_h$.
Because $\nabla\cdot\Uh\subset Q_h$,
$b(\wh,P_h\pbar-\pbar)=0$, and hence
$b(\wh,\eta_p)=b(\wh,\widehat p_h-\pbar)$. The preceding Stokes error
identity implies, for every $\wh\in\Uh$,
\[
  b(\wh,\widehat p_h-\pbar)
  =\gamma a_h(\widehat{\bm v}_h-\vbar,\wh)
   -\gamma\delta_h(\vbar,\wh).
\]
The discrete inf--sup condition, boundedness of $a_h$, and
\eqref{eq:consistency-bound-en} now imply
\begin{align*}
  \norm[0]{\eta_p}
  &\leq
  C\sup_{\bm0\neq\wh\in\Uh}
  \frac{|b(\wh,\eta_p)|}{\tnorm{\wh}}
 \leq
  C\gamma
  \left(
  \tnorm{\widehat{\bm v}_h-\vbar}
  +
  h|\vbar|_2
  \right)
  \leq
  C\gamma h|\vbar|_2.
\end{align*}
Combining this estimate with
$\norm[0]{\pbar-P_h\pbar}\leq Ch|\pbar|_1$
proves \eqref{eq:aux-pressure-error-en}.
\end{proof}

For any $\vh\in\Zh$, we set $T_h=\mathcal S_h(\vh)$ and $q_h=q_h(\vh)$.
We decompose the scalar errors as
\begin{align*}
  &\Tbar-T_h
  =
  (\Tbar-\widehat T_h)
  +
  (\widehat T_h-T_h),\\
  &\qbar-q_h
  =
  (\qbar-\widehat q_h)
  +
  (\widehat q_h-q_h).
\end{align*}
Subtracting the discrete state equation from \eqref{eq:aux-state-en}, adding and subtracting
$c(\vbar;T_h,\phi_h)$, and collecting the unknown difference on the left
leads to
\begin{align}
\label{eq:state-propagation-en}
  &\kappa
  \ip{\nabla(\widehat T_h-T_h)}{\nabla\phi_h}
  +
  c(\vbar;\widehat T_h-T_h,\phi_h)
  =
  -c(\vbar-\vh;T_h,\phi_h)
  \quad
  \forall\phi_h\in S_h,\\
\label{eq:adjoint-propagation-en}
  &\begin{aligned}
  &\kappa
  \ip{\nabla(\widehat q_h-q_h)}{\nabla\psi_h}
  -
  c(\vbar;\widehat q_h-q_h,\psi_h)
  \\
  &\qquad =
  \ip{D(\Tbar-T_h)}{\psi_h}
  +
  c(\vbar-\vh;q_h,\psi_h)
  \quad
  \forall\psi_h\in S_h,
  \end{aligned}
\end{align}

\begin{lemma}
\label{lem:coupling-propagation-en}
For any $\vh\in\Zh$ and its state--adjoint pair $(T_h,q_h)$,
the equations
\eqref{eq:state-propagation-en}--\eqref{eq:adjoint-propagation-en} imply
\begin{align}
 & \norm[1]{\widehat T_h-T_h}
\leq
  C\kappa^{-2}\norm[0]{f}
  \tnorm{\vbar-\vh},
  \label{eq:state-propagation-bound-en}\\
 & \norm[1]{\widehat q_h-q_h}
\leq
  C\kappa^{-1}\norm[1]{\Tbar-T_h}
  +
  C\kappa^{-3}\norm[0]{f}
  \tnorm{\vbar-\vh}.
  \label{eq:adjoint-propagation-bound-en}
\end{align}
\end{lemma}
\begin{proof}
We take
$\phi_h=\widehat T_h-T_h$ in
\eqref{eq:state-propagation-en}. The convection term on the left vanishes,
and Lemma~\ref{lem:discrete-trilinear-en} implies
\begin{align*}
  \kappa\norm[1]{\widehat T_h-T_h}^2
  &\leq
  C\tnorm{\vbar-\vh}
  \norm[1]{T_h}
  \norm[1]{\widehat T_h-T_h}.
\end{align*}
Using \eqref{eq:discrete-state-stability-en} and canceling the error norm
  proves \eqref{eq:state-propagation-bound-en}.

Similarly, we take
$\psi_h=\widehat q_h-q_h$ in
\eqref{eq:adjoint-propagation-en}. Using the $L^2$ contraction of $D$,
Poincar\'e inequality, and the trilinear bound yields
\begin{align*}
  \kappa\norm[1]{\widehat q_h-q_h}^2
  \leq{}&
  C\norm[1]{\Tbar-T_h}
  \norm[1]{\widehat q_h-q_h}
  +
  C\tnorm{\vbar-\vh}
  \norm[1]{q_h}
  \norm[1]{\widehat q_h-q_h}.
\end{align*}
The adjoint bound in \eqref{eq:discrete-adjoint-stability-en} and cancellation,
with the zero-error case immediate, prove
\eqref{eq:adjoint-propagation-bound-en}.
\end{proof}

\subsection{Local error estimates}

We fix a continuous stationary solution satisfying
Assumption~\ref{ass:ssc-en} and use the regularity in
Lemma~\ref{ass:regularity-en}. Since global convexity is not assumed, an arbitrary discrete minimizer need not
approximate the prescribed continuous solution.  For a mesh-independent
radius $r>0$, we select
\begin{equation}
\label{eq:localized-minimizer-en}
 \vbar_h\in\mathop{\arg\min}_{\substack{\wh\in\Zh\\
           \tnorm{\wh-\Pi_h\vbar}\leq r}}j_h(\wh).
\end{equation}
By \eqref{eq:commuting-en},
the interpolated centre $\Pi_h\vbar$ lies in $\Zh$, so the feasible set
is nonempty.
In particular, $\vbar_h-\Pi_h\vbar\in\Zh$ is a valid direction for the
estimates of the first and second derivatives of $j_h$.
Moreover, \eqref{eq:velocity-interpolation-en} implies
$\tnorm{\vbar-\Pi_h\vbar}\leq Ch|\vbar|_2$; hence the triangle inequality
reduces the control error estimate to a bound on
$\tnorm{\vbar_h-\Pi_h\vbar}$.
We therefore first bound the discrete reduced-gradient residual at
$\Pi_h\vbar$ and then establish uniform coercivity of the second
derivative of $j_h$ in this neighbourhood.

We set $\bm v_h^\Pi:=\Pi_h\vbar\in\Zh$ and define
\(T_h^\Pi:=\mathcal S_h(\bm v_h^\Pi),~q_h^\Pi:=q_h(\bm v_h^\Pi).\)
These functions satisfy \eqref{eq:discrete-state-en} and
\eqref{eq:discrete-optimality-adjoint-en} at $\vh=\bm v_h^\Pi$. We write
\begin{equation}
\label{eq:pi-error-notation-en}
 E_T^\Pi:=\norm[1]{\Tbar-T_h^\Pi},
 \qquad E_q^\Pi:=\norm[1]{\qbar-q_h^\Pi}.
\end{equation}

\begin{lemma}[Errors at the interpolated control]
\label{lem:pi-state-adjoint-en}
Under the regularity assumptions of Lemma~\ref{ass:regularity-en},
we have the following estimates
\begin{equation}
\label{eq:pi-state-error-en}
 E_T^\Pi\leq Ch\left(
 M_{\vbar}|\Tbar|_2+\kappa^{-2}\norm[0]{f}|\vbar|_2\right),
\end{equation}
and
\begin{equation}
\label{eq:pi-adjoint-error-en}
 E_q^\Pi\leq Ch\left(
 M_{\vbar}|\qbar|_2+\kappa^{-1}M_{\vbar}|\Tbar|_2
 +\kappa^{-3}\norm[0]{f}|\vbar|_2\right).
\end{equation}
\end{lemma}
\begin{proof}
Using the auxiliary variables in
\eqref{eq:aux-state-en}--\eqref{eq:aux-adjoint-en}, we decompose the errors as
\begin{align*}
 \Tbar-T_h^\Pi
 &=(\Tbar-\widehat T_h)+(\widehat T_h-T_h^\Pi),\\
 \qbar-q_h^\Pi
 &=(\qbar-\widehat q_h)+(\widehat q_h-q_h^\Pi).
\end{align*}
Lemma~\ref{lem:auxiliary-errors-en} bounds the first term in each
decomposition.  For the second terms, we specialize the comparison
equations underlying Lemma~\ref{lem:coupling-propagation-en} to
$\vh=\bm v_h^\Pi$.
Since $\bm v_h^\Pi\in\Zh$, the stability bounds
\eqref{eq:discrete-state-stability-en} and
\eqref{eq:discrete-adjoint-stability-en} apply to
$T_h^\Pi$ and $q_h^\Pi$, respectively.
Moreover, \eqref{eq:velocity-interpolation-en} implies
$\tnorm{\vbar-\bm v_h^\Pi}\leq Ch|\vbar|_2$.

We subtract the state equation at $\bm v_h^\Pi$ from
\eqref{eq:aux-state-en} and collect the terms with velocity $\vbar$
on the left to obtain
\[
 \begin{aligned}
 &\kappa\ip{\nabla(\widehat T_h-T_h^\Pi)}{\nabla\phi_h}
       +c(\vbar;\widehat T_h-T_h^\Pi,\phi_h)=-c(\vbar-\bm v_h^\Pi;T_h^\Pi,\phi_h)
       \qquad\forall\phi_h\in S_h.
 \end{aligned}
\]
We choose $\phi_h=\widehat T_h-T_h^\Pi$.  The convection term on the left
vanishes because $\vbar\in\Z$, and Lemma~\ref{lem:discrete-trilinear-en}
therefore gives
\[
 \kappa\norm[1]{\widehat T_h-T_h^\Pi}^2
 \leq C\tnorm{\vbar-\bm v_h^\Pi}\norm[1]{T_h^\Pi}
             \norm[1]{\widehat T_h-T_h^\Pi}.
\]
Then we further derive
\begin{align*}
 \norm[1]{\widehat T_h-T_h^\Pi}
 &\leq C\kappa^{-1}\tnorm{\vbar-\bm v_h^\Pi}\norm[1]{T_h^\Pi}\\
 &\leq C\kappa^{-2}\norm[0]{f}\tnorm{\vbar-\bm v_h^\Pi}
 \leq Ch\kappa^{-2}\norm[0]{f}|\vbar|_2.
\end{align*}
The triangle inequality and \eqref{eq:aux-state-error-en} now imply
\begin{align*}
 E_T^\Pi
 &\leq \norm[1]{\Tbar-\widehat T_h}
       +\norm[1]{\widehat T_h-T_h^\Pi}\\
 &\leq Ch\left(M_{\vbar}|\Tbar|_2
       +\kappa^{-2}\norm[0]{f}|\vbar|_2\right).
\end{align*}
This proves \eqref{eq:pi-state-error-en}.

For the adjoint, we subtract
\eqref{eq:discrete-optimality-adjoint-en} at $\bm v_h^\Pi$ from
\eqref{eq:aux-adjoint-en}.  The difference satisfies
\[
 \begin{aligned}
 &\kappa\ip{\nabla(\widehat q_h-q_h^\Pi)}{\nabla\psi_h}
       -c(\vbar;\widehat q_h-q_h^\Pi,\psi_h)\\
  &\qquad =\ip{D(\Tbar-T_h^\Pi)}{\psi_h}
       +c(\vbar-\bm v_h^\Pi;q_h^\Pi,\psi_h)
       \quad\forall\psi_h\in S_h.
 \end{aligned}
\]
Testing with $\psi_h=\widehat q_h-q_h^\Pi$ again removes the convection
term on the left.  The contraction of $D$, Poincar\'e  inequality, and
Lemma~\ref{lem:discrete-trilinear-en} imply
\begin{align*}
 \kappa\norm[1]{\widehat q_h-q_h^\Pi}^2
 \leq{}&C E_T^\Pi\norm[1]{\widehat q_h-q_h^\Pi}+C\tnorm{\vbar-\bm v_h^\Pi}\norm[1]{q_h^\Pi}
              \norm[1]{\widehat q_h-q_h^\Pi}.
\end{align*}
The adjoint stability bound gives
\begin{align*}
 \norm[1]{\widehat q_h-q_h^\Pi}
 &\leq C\kappa^{-1}E_T^\Pi
       +C\kappa^{-1}\tnorm{\vbar-\bm v_h^\Pi}\norm[1]{q_h^\Pi}\\
 &\leq C\kappa^{-1}E_T^\Pi
       +C\kappa^{-3}\norm[0]{f}\tnorm{\vbar-\bm v_h^\Pi}.
\end{align*}
Combining this inequality with \eqref{eq:aux-adjoint-error-en} and the
velocity interpolation bound, we obtain
\begin{align*}
 E_q^\Pi
 &\leq \norm[1]{\qbar-\widehat q_h}
       +\norm[1]{\widehat q_h-q_h^\Pi}\\
 &\leq ChM_{\vbar}|\qbar|_2+C\kappa^{-1}E_T^\Pi
       +Ch\kappa^{-3}\norm[0]{f}|\vbar|_2\\
 &\leq ChM_{\vbar}|\qbar|_2
       +Ch\kappa^{-1}\left(M_{\vbar}|\Tbar|_2
                    +\kappa^{-2}\norm[0]{f}|\vbar|_2\right)+Ch\kappa^{-3}\norm[0]{f}|\vbar|_2\\
 &\leq Ch\left(M_{\vbar}|\qbar|_2
       +\kappa^{-1}M_{\vbar}|\Tbar|_2
       +\kappa^{-3}\norm[0]{f}|\vbar|_2\right).
\end{align*}
Here the third inequality uses \eqref{eq:pi-state-error-en}; the last
combines the two contributions proportional to
$h\kappa^{-3}\norm[0]{f}|\vbar|_2$ into the generic constant $C$.
This proves \eqref{eq:pi-adjoint-error-en}.
\end{proof}

We define the data and regularity factor
\begin{equation}
\label{eq:residual-data-factor-en}
\begin{aligned}
  \mathcal M_{\mathrm{res}}
  :={}&
  \left(
  \gamma+\kappa^{-4}\norm[0]{f}^2
  \right)|\vbar|_2+
  \kappa^{-2}M_{\vbar}\norm[0]{f}|\Tbar|_2
  +
  \kappa^{-1}M_{\vbar}\norm[0]{f}|\qbar|_2.
\end{aligned}
\end{equation}

For $\ell_h\in\Zh'$, we set
\begin{equation}
\label{eq:discrete-dual-norm-en}
 \dualnorm{\ell_h}:=
 \sup_{\bm0\ne\zh\in\Zh}\frac{|\ell_h(\zh)|}{\tnorm{\zh}}.
\end{equation}

\begin{lemma}[Discrete reduced-gradient residual]
\label{lem:discrete-gradient-residual-en}
The residual of the discrete reduced gradient at the interpolated continuous
control satisfies
\begin{equation}
\label{eq:gradient-residual-bound-en}
  \dualnorm{j_h'(\Pi_h\vbar)}
  \leq
  Ch\mathcal M_{\mathrm{res}}.
\end{equation}
\end{lemma}

\begin{proof}
We fix $\zh\in\Zh$.  Proposition~\ref{prop:discrete-mass-conservation-en}
implies $\nabla\cdot\zh=0$, and hence
\[
 b(\zh,\pbar)=\ip{\nabla\cdot\zh}{\pbar}=0.
\]
The continuous control identity \eqref{eq:continuous-control-tested-en}
therefore reduces to
\[
 \gamma a_h(\vbar,\zh)-\ip{\qbar\nabla\Tbar}{\zh}
 =-\gamma\delta_h(\vbar,\zh).
\]
Applying \eqref{eq:discrete-gradient-en} at $\bm v_h^\Pi=\Pi_h\vbar$
and adding and subtracting the continuous control terms, we obtain
\[
 \begin{aligned}
 j_h'(\Pi_h\vbar)\zh
 ={}&\gamma a_h(\Pi_h\vbar-\vbar,\zh)+\left[\gamma a_h(\vbar,\zh)
         -\ip{\qbar\nabla\Tbar}{\zh}\right]\\
    &+\ip{\qbar\nabla\Tbar-q_h^\Pi\nabla T_h^\Pi}{\zh}.
 \end{aligned}
\]
Substitution of the preceding continuous identity establishes
\begin{align}
 j_h'(\Pi_h\vbar)\zh
 ={}&\gamma a_h(\Pi_h\vbar-\vbar,\zh)
      -\gamma\delta_h(\vbar,\zh)+\ip{\qbar\nabla\Tbar-q_h^\Pi\nabla T_h^\Pi}{\zh}.
 \label{eq:gradient-residual-decomposition-en}
\end{align}
The first two terms are the interpolation and consistency contributions.
By the Cauchy--Schwarz inequality for $a_h$ and
\eqref{eq:velocity-interpolation-en}, we obtain
\begin{align*}
 \gamma|a_h(\Pi_h\vbar-\vbar,\zh)|
 &\leq\gamma\tnorm{\Pi_h\vbar-\vbar}\tnorm{\zh}\leq C\gamma h|\vbar|_2\tnorm{\zh},
\end{align*}
while \eqref{eq:consistency-bound-en} implies
\[
 \gamma|\delta_h(\vbar,\zh)|
 \leq C\gamma h|\vbar|_2\tnorm{\zh}.
\]

For the nonlinear contribution,
H\"older's inequality, the embedding
$H_0^1(\Omega)\hookrightarrow L^4(\Omega)$, and
Lemma~\ref{lem:discrete-trilinear-en} imply
\begin{align*}
 \left|\ip{\qbar\nabla\Tbar-q_h^\Pi\nabla T_h^\Pi}{\zh}\right|
 \leq{}&\norm[L^4(\Omega)]{\qbar-q_h^\Pi}\norm[1]{\Tbar}
          \norm[L^4(\Omega)]{\zh}\\*
       &+\norm[L^4(\Omega)]{q_h^\Pi}\norm[1]{\Tbar-T_h^\Pi}
          \norm[L^4(\Omega)]{\zh}\\*
 \leq{}&C\left(E_q^\Pi\norm[1]{\Tbar}
       +\norm[1]{q_h^\Pi}E_T^\Pi\right)\tnorm{\zh}.
\end{align*}
By \eqref{eq:state-stability-en} and
\eqref{eq:discrete-adjoint-stability-en}, we obtain
\[
 E_q^\Pi\norm[1]{\Tbar}+\norm[1]{q_h^\Pi}E_T^\Pi
 \leq C\norm[0]{f}
       \left(\kappa^{-1}E_q^\Pi+\kappa^{-2}E_T^\Pi\right).
\]
We estimate these two errors separately. From
\eqref{eq:pi-adjoint-error-en}, we obtain
\begin{align*}
 \kappa^{-1}\norm[0]{f}E_q^\Pi
 &\leq Ch\kappa^{-1}\norm[0]{f}
       \left(M_{\vbar}|\qbar|_2
       +\kappa^{-1}M_{\vbar}|\Tbar|_2
       +\kappa^{-3}\norm[0]{f}|\vbar|_2\right)\\*
 &\leq Ch\Big[
       \kappa^{-1}M_{\vbar}\norm[0]{f}|\qbar|_2
       +\kappa^{-2}M_{\vbar}\norm[0]{f}|\Tbar|_2
       +\kappa^{-4}\norm[0]{f}^2|\vbar|_2\Big].
\end{align*}
Similarly, \eqref{eq:pi-state-error-en} implies
\begin{align*}
 \kappa^{-2}\norm[0]{f}E_T^\Pi
 &\leq Ch\kappa^{-2}\norm[0]{f}
       \left(M_{\vbar}|\Tbar|_2
             +\kappa^{-2}\norm[0]{f}|\vbar|_2\right)\\*
 &\leq Ch\Big[
       \kappa^{-2}M_{\vbar}\norm[0]{f}|\Tbar|_2
       +\kappa^{-4}\norm[0]{f}^2|\vbar|_2\Big].
\end{align*}
The terms proportional to $\kappa^{-2}M_{\vbar}\norm[0]{f}|\Tbar|_2$
and $\kappa^{-4}\norm[0]{f}^2|\vbar|_2$ occur in both bounds;
their multiplicities are absorbed into $C$.  It follows that
\begin{align*}
 &\left|\ip{\qbar\nabla\Tbar-q_h^\Pi\nabla T_h^\Pi}{\zh}\right|\\*
 &\qquad\leq Ch\Big[
       \kappa^{-1}M_{\vbar}\norm[0]{f}|\qbar|_2
       +\kappa^{-2}M_{\vbar}\norm[0]{f}|\Tbar|_2
       +\kappa^{-4}\norm[0]{f}^2|\vbar|_2\Big]\tnorm{\zh}.
\end{align*}

Combining the three contributions in
\eqref{eq:gradient-residual-decomposition-en}, we obtain
\begin{align*}
 |j_h'(\Pi_h\vbar)\zh|
 &\leq Ch\Big[
       \bigl(\gamma+\kappa^{-4}\norm[0]{f}^2\bigr)|\vbar|_2
       +\kappa^{-2}M_{\vbar}\norm[0]{f}|\Tbar|_2\\*
 &\qquad\quad
       +\kappa^{-1}M_{\vbar}\norm[0]{f}|\qbar|_2\Big]\tnorm{\zh}\\*
 &=Ch\mathcal M_{\mathrm{res}}\tnorm{\zh},
\end{align*}
where the last equality uses \eqref{eq:residual-data-factor-en}.
For $\zh\ne\bm0$, dividing by $\tnorm{\zh}$ and taking the supremum
in \eqref{eq:discrete-dual-norm-en} proves
\eqref{eq:gradient-residual-bound-en}.  The constant is independent of
$h$, and no pressure approximation term enters the estimate because
$b(\zh,\pbar)$ vanishes exactly.
\end{proof}

\begin{lemma}
\label{lem:discrete-second-derivative-lipschitz-en}
There is a constant $L=C\kappa^{-5}\norm[0]{f}^2$, independent of $h$, such that
\begin{equation}
\label{eq:discrete-second-derivative-lipschitz-en}
 \left|(j_h''(\bm v_{1,h})-j_h''(\bm v_{2,h}))[\bm d_h,\bm d_h]\right|
 \leq L\tnorm{\bm v_{1,h}-\bm v_{2,h}}\tnorm{\bm d_h}^2
\end{equation}
for all $\bm v_{1,h},\bm v_{2,h},\bm d_h\in\Zh$.
\end{lemma}
\begin{proof}
We write $T_i=\mathcal S_h(\bm v_{i,h})$,
$q_i=q_h(\bm v_{i,h})$, and
$z_i=\mathcal S_h'(\bm v_{i,h})\bm d_h$, $i=1,2$.
The estimates \eqref{eq:discrete-state-stability-en},
\eqref{eq:discrete-adjoint-stability-en}, and
\eqref{eq:discrete-linearized-state-bound-en} imply
\[
 \norm[1]{T_i}\leq C\kappa^{-1}\norm[0]{f},\qquad
 \norm[1]{q_i}\leq C\kappa^{-2}\norm[0]{f},\qquad
 \norm[1]{z_i}\leq C\kappa^{-2}\norm[0]{f}\tnorm{\bm d_h}.
\]
For $A_h$ defined in \eqref{eq:discrete-convection-diffusion-form-en},
skew-symmetry implies
$A_h(\bm v_{1,h};\phi_h,\phi_h)=\kappa\norm[1]{\phi_h}^2$.
Moreover, the bilinear form satisfies
\[
 A_h(\bm v_{1,h};\eta_h,\phi_h)
 -A_h(\bm v_{2,h};\eta_h,\phi_h)
 =c(\bm v_{1,h}-\bm v_{2,h};\eta_h,\phi_h).
\]
Taking differences in \eqref{eq:discrete-state-en},
\eqref{eq:discrete-optimality-adjoint-en}, and
\eqref{eq:discrete-linearized-state-en} at $\bm v_{1,h}$ and $\bm v_{2,h}$,
we obtain the following identities for every $\phi_h\in S_h$:
\begin{align*}
 A_h(\bm v_{1,h};T_1-T_2,\phi_h)
 &=-c(\bm v_{1,h}-\bm v_{2,h};T_2,\phi_h),\\
 A_h(\bm v_{1,h};\phi_h,q_1-q_2)
 &=\ip{D(T_1-T_2)}{\phi_h}
   +c(\bm v_{1,h}-\bm v_{2,h};q_2,\phi_h),\\*
 A_h(\bm v_{1,h};z_1-z_2,\phi_h)
 &=-c(\bm v_{1,h}-\bm v_{2,h};z_2,\phi_h)-c(\bm d_h;T_1-T_2,\phi_h).
\end{align*}
In the adjoint identity, we used
$A_h(\bm v_{i,h};\phi_h,q_i)=\ip{DT_i}{\phi_h}$ and
\[
 -c(\bm v_{1,h}-\bm v_{2,h};\phi_h,q_2)
 =c(\bm v_{1,h}-\bm v_{2,h};q_2,\phi_h),
\]
which follows from $\bm v_{1,h}-\bm v_{2,h}\in\Zh$
and the skew-symmetry of $c$.

Testing the state-difference equation with $T_1-T_2$ and applying
Lemma~\ref{lem:discrete-trilinear-en}, we obtain
\[
 \kappa\norm[1]{T_1-T_2}^2
 \leq C\tnorm{\bm v_{1,h}-\bm v_{2,h}}
          \norm[1]{T_2}\norm[1]{T_1-T_2}.
\]
After cancelling the nonzero error norm and applying the state stability bound, we obtain
\begin{align*}
 \norm[1]{T_1-T_2}
 &\leq C\kappa^{-1}\tnorm{\bm v_{1,h}-\bm v_{2,h}}\norm[1]{T_2}\leq C\kappa^{-2}\norm[0]{f}\tnorm{\bm v_{1,h}-\bm v_{2,h}}.
\end{align*}

Next, we test the adjoint-difference equation with $q_1-q_2$.
The contraction property of $D$, Poincar\'e  inequality, and
Lemma~\ref{lem:discrete-trilinear-en} imply
\[
 \begin{aligned}
 \kappa\norm[1]{q_1-q_2}^2
 &\leq\norm[0]{D(T_1-T_2)}\norm[0]{q_1-q_2}+C\tnorm{\bm v_{1,h}-\bm v_{2,h}}
             \norm[1]{q_2}\norm[1]{q_1-q_2}\\
 &\leq C\Big(\norm[1]{T_1-T_2}
       +\tnorm{\bm v_{1,h}-\bm v_{2,h}}\norm[1]{q_2}\Big)
       \norm[1]{q_1-q_2}.
 \end{aligned}
\]
Dividing by $\kappa\norm[1]{q_1-q_2}$ and substituting the state-difference
and adjoint stability bounds yields
\begin{align*}
 \norm[1]{q_1-q_2}
 &\leq C\kappa^{-1}\Big(\norm[1]{T_1-T_2}
       +\tnorm{\bm v_{1,h}-\bm v_{2,h}}\norm[1]{q_2}\Big)\\*
 &\leq C\kappa^{-3}\norm[0]{f}\tnorm{\bm v_{1,h}-\bm v_{2,h}}.
\end{align*}

For the linearized-state difference, we test with $z_1-z_2$ to obtain
\[
 \kappa\norm[1]{z_1-z_2}^2
 \leq C\Big(\tnorm{\bm v_{1,h}-\bm v_{2,h}}\norm[1]{z_2}
        +\tnorm{\bm d_h}\norm[1]{T_1-T_2}\Big)\norm[1]{z_1-z_2}.
\]
Consequently,
\begin{align*}
 \norm[1]{z_1-z_2}
 &\leq C\kappa^{-1}\Big(
      \tnorm{\bm v_{1,h}-\bm v_{2,h}}\norm[1]{z_2}
      +\tnorm{\bm d_h}\norm[1]{T_1-T_2}\Big)\\*
 &\leq C\kappa^{-3}\norm[0]{f}
       \tnorm{\bm v_{1,h}-\bm v_{2,h}}\tnorm{\bm d_h}.
\end{align*}
The last inequality uses the linearized-state stability bound for $z_2$
and the state-difference estimate; both contributions have the same
dependence on $\kappa$.

It remains to compare the second derivatives of $j_h$ at the two controls.
Their identical regularization terms cancel, and bilinearity implies
\[
 \begin{aligned}
 &(j_h''(\bm v_{1,h})-j_h''(\bm v_{2,h}))[\bm d_h,\bm d_h]\\
 &\qquad=\norm[0]{Dz_1}^2-\norm[0]{Dz_2}^2
       +2c(\bm d_h;q_1,z_1)-2c(\bm d_h;q_2,z_2)\\
 &\qquad=\ip{D(z_1+z_2)}{D(z_1-z_2)}
       +2c(\bm d_h;q_1,z_1-z_2)
       +2c(\bm d_h;q_1-q_2,z_2).
 \end{aligned}
\]
For the first term, the contraction of $D$ and Poincar\'e inequality
reduce the estimate to the bounds for $z_i$ and $z_1-z_2$:
\begin{align*}
 \left|\ip{D(z_1+z_2)}{D(z_1-z_2)}\right|
 &\leq C\big(\norm[1]{z_1}+\norm[1]{z_2}\big)\norm[1]{z_1-z_2}\\*
 &\leq C\kappa^{-5}\norm[0]{f}^2
       \tnorm{\bm v_{1,h}-\bm v_{2,h}}\tnorm{\bm d_h}^2.
\end{align*}
The two convection terms are estimated separately by
Lemma~\ref{lem:discrete-trilinear-en}:
\begin{align*}
 2|c(\bm d_h;q_1,z_1-z_2)|
 &\leq C\tnorm{\bm d_h}\norm[1]{q_1}\norm[1]{z_1-z_2}\\*
 &\leq C\kappa^{-5}\norm[0]{f}^2
       \tnorm{\bm v_{1,h}-\bm v_{2,h}}\tnorm{\bm d_h}^2,\\
 2|c(\bm d_h;q_1-q_2,z_2)|
 &\leq C\tnorm{\bm d_h}\norm[1]{q_1-q_2}\norm[1]{z_2}\\*
 &\leq C\kappa^{-5}\norm[0]{f}^2
       \tnorm{\bm v_{1,h}-\bm v_{2,h}}\tnorm{\bm d_h}^2.
\end{align*}
Adding these three estimates proves
\eqref{eq:discrete-second-derivative-lipschitz-en} with
$L=C\kappa^{-5}\norm[0]{f}^2$.
Here $C$ depends only on the domain, mesh shape regularity, and the
uniform spectral-equivalence constants of $a_R$; it is independent of
$h$, $\kappa$, $\gamma$, $f$, the two control points, and $\bm d_h$.
\end{proof}
\begin{lemma}
\label{lem:discrete-direction-compactness-en}
Let $h\to0$ along the uniformly shape-regular mesh family.
If $\bm w_h=\bm w_h^1+\bm w_h^R\in\Zh$ is bounded in $\tnorm{\cdot}$,
then, along a subsequence,
\[
 \bm w_h^1\rightharpoonup\bm w\text{ in }[H_0^1(\Omega)]^d,
 \qquad \bm w_h\to\bm w\text{ in }[L^4(\Omega)]^d,
 \qquad \bm w\in\Z,
\]
and $\norm[1]{\bm w}\leq\liminf_{h\to0}\tnorm{\bm w_h}$.
\end{lemma}
\begin{proof}
We first select a subsequence along which
$\tnorm{\bm w_h}$ converges to its lower limit.
Since $h_K\leq h$, the lower bound in
\eqref{eq:aR-spectral-equivalence-en} and the definitions
\eqref{eq:ah-en} and \eqref{eq:velocity-energy-norm-en} imply
\begin{align*}
 \norm[0]{\bm w_h^R}^2
 &=\sum_{K\in\Th}h_K^2
       \bigl(h_K^{-2}\norm[0,K]{\bm w_h^R}^2\bigr)\\
 &\leq h^2\sum_{K\in\Th}h_K^{-2}\norm[0,K]{\bm w_h^R}^2
 \leq c_R^{-1}h^2 a_R(\bm w_h^R,\bm w_h^R)\\
 &\leq c_R^{-1}h^2
       \bigl(\norm[1]{\bm w_h^1}^2+a_R(\bm w_h^R,\bm w_h^R)\bigr)
 =c_R^{-1}h^2\tnorm{\bm w_h}^2.
\end{align*}
Taking square roots proves
\[
 \norm[0]{\bm w_h^R}\leq c_R^{-1/2}h\tnorm{\bm w_h},
\]
where $c_R$ is independent of $h$.
The inverse estimate in Lemma~\ref{lem:discrete-trilinear-en}, together
with \eqref{eq:aR-spectral-equivalence-en}, also implies
\[
 \norm[L^4(\Omega)]{\bm w_h^R}
 \leq Ch^{1-d/4}\tnorm{\bm w_h}\longrightarrow0.
\]
The conforming components are bounded in $[H_0^1(\Omega)]^d$.
Weak compactness and $H_0^1(\Omega)\Subset L^4(\Omega)$ for $d\leq3$
therefore establish the stated convergences.  Passing to the limit in
$\ip{\bm w_h}{\nabla\phi}=0$, $\phi\in C_c^\infty(\Omega)$,
shows that $\bm w$ is divergence-free.  Its boundary trace is zero by weak
$H_0^1$ convergence.  The energy bound follows from weak lower semicontinuity.
\end{proof}
\begin{lemma}[Local discrete coercivity]
\label{lem:directional-coercivity-en}
Under Assumption~\ref{ass:ssc-en} and the regularity setting
of Lemma~\ref{ass:regularity-en}, there exist $r>0$, $h_0>0$, and
$\mu>0$, independent of $h$, such that
\begin{equation}
\label{eq:directional-coercivity-en}
 j_h''(\wh)[\bm d_h,\bm d_h]\geq\mu\tnorm{\bm d_h}^2
 \quad\forall\wh,\bm d_h\in\Zh,
 \quad\tnorm{\wh-\Pi_h\vbar}\leq r,
 \quad h<h_0.
\end{equation}
One may take $\mu=\frac14\min\{\delta,\gamma\}$.
\end{lemma}
\begin{proof}
We set $m=\min\{\delta,\gamma\}$.  We first prove coercivity at
$\bm v_h^\Pi=\Pi_h\vbar$ with constant $m/2$.
If this assertion were false, there would exist a sequence of mesh sizes
$h\to0$ and corresponding directions $\bm d_h\in\Zh$ with
$\tnorm{\bm d_h}=1$ and
$j_h''(\bm v_h^\Pi)[\bm d_h,\bm d_h]<m/2$.
Lemma~\ref{lem:discrete-direction-compactness-en} supplies a subsequence
such that $\bm d_h\to\bm d$ in $L^4$, where $\bm d\in\Z$ and
$\norm[1]{\bm d}\leq1$.
Lemma~\ref{lem:pi-state-adjoint-en} implies
$T_h^\Pi\to\Tbar$ and $q_h^\Pi\to\qbar$ in $H_0^1$.

We set $z_h=\mathcal S_h'(\bm v_h^\Pi)\bm d_h$ and
$z=S'(\vbar)\bm d$.
The density of $C_0^\infty(\Omega)$ in $H_0^1(\Omega)$ and
\eqref{eq:scalar-interpolation-en} imply
$\inf_{\phi_h\in S_h}\norm[1]{z-\phi_h}\to0$.
We choose $\zeta_h\in S_h$ such that $\norm[1]{z-\zeta_h}\to0$.
Moreover, \eqref{eq:velocity-interpolation-en} implies
\[
 \tnorm{\bm v_h^\Pi}
 \leq\norm[1]{\vbar}+Ch|\vbar|_2,
 \qquad\tnorm{\bm v_h^\Pi-\vbar}\longrightarrow0.
\]
Thus the continuity constants of $A_h(\bm v_h^\Pi;\cdot,\cdot)$
are bounded independently of $h$.
Subtracting the two linearized equations gives
\begin{align*}
 &A_h(\bm v_h^\Pi;z_h-\zeta_h,\phi_h)\\*
 &\qquad=A_h(\bm v_h^\Pi;z-\zeta_h,\phi_h)
       -c(\bm v_h^\Pi-\vbar;z,\phi_h)\\*
 &\qquad\quad -c(\bm d_h-\bm d;\Tbar,\phi_h)
       -c(\bm d_h;T_h^\Pi-\Tbar,\phi_h).
\end{align*}
Testing with $\phi_h=z_h-\zeta_h$, and using coercivity, H\"older's
inequality, we obtain
\begin{align}
 \kappa\norm[1]{z_h-\zeta_h}
 \leq{}&\bigl(\kappa+C\tnorm{\bm v_h^\Pi}\bigr)
             \norm[1]{z-\zeta_h}\notag+C\tnorm{\bm v_h^\Pi-\vbar}\norm[1]{z}\notag\\*
 &+C\norm[L^4(\Omega)]{\bm d_h-\bm d}\norm[1]{\Tbar}\notag+C\tnorm{\bm d_h}\norm[1]{T_h^\Pi-\Tbar}.
 \label{eq:linearized-state-strong-convergence-en}
\end{align}
Each term on the right tends to zero: the first two by approximation,
the third by Lemma~\ref{lem:discrete-direction-compactness-en}, and
the fourth by Lemma~\ref{lem:pi-state-adjoint-en} and
$\tnorm{\bm d_h}=1$.  Consequently,
\[
 \norm[1]{z_h-z}
 \leq\norm[1]{z_h-\zeta_h}+\norm[1]{\zeta_h-z}
 \longrightarrow0.
\]
To pass to the limit in the second derivative of $j_h$, we recall from
\eqref{eq:discrete-second-derivative-en} that
\[
 j_h''(\bm v_h^\Pi)[\bm d_h,\bm d_h]
 =\norm[0]{Dz_h}^2
  +2\ip{z_h}{\bm d_h\cdot\nabla q_h^\Pi}
  +\gamma\tnorm{\bm d_h}^2.
\]
In particular, $z_h\to z$ in $L^4(\Omega)$, and
\begin{align*}
 \norm[0]{z_h\bm d_h-z\bm d}
 &\leq\norm[L^4(\Omega)]{z_h-z}\norm[L^4(\Omega)]{\bm d_h}+\norm[L^4(\Omega)]{z}\norm[L^4(\Omega)]{\bm d_h-\bm d}
 \longrightarrow0.
\end{align*}
Together with $q_h^\Pi\to\qbar$ in $H_0^1(\Omega)$, this yields
\begin{align*}
 &\left|\ip{z_h}{\bm d_h\cdot\nabla q_h^\Pi}
              -\ip z{\bm d\cdot\nabla\qbar}\right|\leq\norm[0]{z_h\bm d_h-z\bm d}\norm[1]{q_h^\Pi}
       +\norm[0]{z\bm d}\norm[1]{q_h^\Pi-\qbar}
 \longrightarrow0.
\end{align*}
Also, $Dz_h\to Dz$ in $L^2(\Omega)$.
Consequently, \eqref{eq:discrete-second-derivative-en} implies
\begin{align*}
 \lim_{h\to0}j_h''(\bm v_h^\Pi)[\bm d_h,\bm d_h]
 &=\norm[0]{Dz}^2+2\ip{z}{\bm d\cdot\nabla\qbar}+\gamma\\
 &=j''(\vbar)[\bm d,\bm d]
       +\gamma\bigl(1-\norm[1]{\bm d}^2\bigr)\\
 &\geq\delta\norm[1]{\bm d}^2
       +\gamma\bigl(1-\norm[1]{\bm d}^2\bigr)\geq m.
\end{align*}
The last expression is a convex combination of $\delta$ and $\gamma$.
Together with the assumed upper bound, this implies
\[
 m\leq\lim_{h\to0}
 j_h''(\bm v_h^\Pi)[\bm d_h,\bm d_h]
 \leq\frac m2.
\]
This contradiction establishes coercivity at the centre for all
sufficiently small $h$.
We choose $r>0$ with $Lr\leq m/4$, where $L$ is from
Lemma~\ref{lem:discrete-second-derivative-lipschitz-en}. %
For $\wh\in\Zh$ with $\tnorm{\wh-\Pi_h\vbar}\leq r$
and every $\bm d_h\in\Zh$, that lemma implies
\begin{align*}
 j_h''(\wh)[\bm d_h,\bm d_h]
 &\geq j_h''(\Pi_h\vbar)[\bm d_h,\bm d_h]
       -L\tnorm{\wh-\Pi_h\vbar}\tnorm{\bm d_h}^2\\*
 &\geq\left(\frac m2-Lr\right)\tnorm{\bm d_h}^2
 \geq\frac m4\tnorm{\bm d_h}^2.
\end{align*}
This proves \eqref{eq:directional-coercivity-en} with $\mu=m/4$.
\end{proof}

We use the following error notation:
\begin{equation}
\label{eq:energy-error-notation-en}
 E_v=\tnorm{\vbar-\vbar_h},\quad
 E_T=\norm[1]{\Tbar-T_h},\quad
 E_q=\norm[1]{\qbar-q_h},\quad
 E_p=\norm[0]{\pbar-p_h}.
\end{equation}

\begin{theorem}[Local existence, uniqueness, and control error]
\label{thm:control-error-en}
Under Assumption~\ref{ass:ssc-en} and the regularity setting of
Lemma~\ref{ass:regularity-en},
 with $r$ chosen as in Lemma~\ref{lem:directional-coercivity-en},
problem \eqref{eq:localized-minimizer-en}
has a unique minimizer for all sufficiently small $h$.
Moreover,
\begin{equation}
\label{eq:gradient-error-relation-en}
 \tnorm{\vbar_h-\Pi_h\vbar}
 \leq\mu^{-1}\dualnorm{j_h'(\Pi_h\vbar)}
 \leq Ch\mu^{-1}\mathcal M_{\mathrm{res}},
\end{equation}
and
\begin{equation}
\label{eq:control-error-final-en}
 E_v\leq Ch|\vbar|_2+Ch\mu^{-1}\mathcal M_{\mathrm{res}}.
\end{equation}
\end{theorem}
\begin{proof}
We take $r,h_0,\mu>0$ as in
Lemma~\ref{lem:directional-coercivity-en}, fix $h<h_0$, and set
\[
 \mathcal B_h:=\{\wh\in\Zh:
                 \tnorm{\wh-\Pi_h\vbar}\leq r\}.
\]
Since $\Pi_h\vbar\in\Zh$ by \eqref{eq:commuting-en}, this ball is
nonempty and compact in the finite-dimensional space $\Zh$.
Continuity of $j_h$ therefore ensures the existence of a minimizer.
By the convexity of $\mathcal B_h$, Taylor's formula and
Lemma~\ref{lem:directional-coercivity-en} imply
\[
 j_h(\wh)\geq j_h(\bm u_h)
 +j_h'(\bm u_h)(\wh-\bm u_h)
 +\frac\mu2\tnorm{\wh-\bm u_h}^2
 \qquad\forall\bm u_h,\wh\in\mathcal B_h.
\]
Thus $j_h$ is strongly convex on $\mathcal B_h$, and its minimizer
$\vbar_h$ is unique.  Since the ball is convex, the constrained
first-order condition reads
\begin{equation}
\label{eq:localized-vi-en}
 j_h'(\vbar_h)(\wh-\vbar_h)\geq0
 \qquad\forall\wh\in\Zh\text{ with }\tnorm{\wh-\Pi_h\vbar}\leq r.
\end{equation}

We set $\bm e_h=\vbar_h-\Pi_h\vbar\in\Zh$.
Taking $\wh=\Pi_h\vbar$ in \eqref{eq:localized-vi-en} yields
$j_h'(\vbar_h)\bm e_h\leq0$.
Since $\mathcal B_h$ is convex,
$\Pi_h\vbar+t\bm e_h\in\mathcal B_h$ for every $t\in[0,1]$.
Hence Lemma~\ref{lem:directional-coercivity-en} applies at
$\Pi_h\vbar+t\bm e_h$ for each such $t$, and the integral mean-value
identity implies
\begin{align*}
 \mu\tnorm{\bm e_h}^2
 &\leq\int_0^1j_h''(\Pi_h\vbar+t\bm e_h)
                 [\bm e_h,\bm e_h]\,\mathrm dt\\
 &=\bigl[j_h'(\vbar_h)-j_h'(\Pi_h\vbar)\bigr]\bm e_h\\
 &\leq-j_h'(\Pi_h\vbar)\bm e_h\\
 &\leq\dualnorm{j_h'(\Pi_h\vbar)}\tnorm{\bm e_h}.
\end{align*}
By
Lemma~\ref{lem:discrete-gradient-residual-en} we can establish
\[
 \tnorm{\bm e_h}
 \leq\mu^{-1}\dualnorm{j_h'(\Pi_h\vbar)}
 \leq Ch\mu^{-1}\mathcal M_{\mathrm{res}}.
\]
The triangle inequality and \eqref{eq:velocity-interpolation-en} then imply
\begin{align*}
 E_v=\tnorm{\vbar-\vbar_h}
 &\leq\tnorm{\vbar-\Pi_h\vbar}+\tnorm{\bm e_h}\\
 &\leq Ch|\vbar|_2+ Ch\mu^{-1} \mathcal M_{\mathrm{res}},
\end{align*}
which proves \eqref{eq:control-error-final-en}.

Finally, $r$ is fixed, whereas
$Ch\mu^{-1} \mathcal M_{\mathrm{res}} \to0$ as $h\to0$.
Thus $\tnorm{\vbar_h-\Pi_h\vbar}<r$ for sufficiently small $h$,
and $\vbar_h$ is interior to $\mathcal B_h$ relative to $\Zh$.
For every $\zh\in\Zh$, the perturbations $\vbar_h\pm t\zh$
are therefore feasible for sufficiently small $t>0$.
Applying \eqref{eq:localized-vi-en} to both signs proves
\[
 j_h'(\vbar_h)\zh=0\qquad\forall\zh\in\Zh.
\]
The strong-convexity inequality above then establishes strict local
minimality on $\Zh$ and excludes any other stationary point in
$\mathcal B_h$.
\end{proof}

By Theorem~\ref{thm:control-error-en}, the selected control is an
unconstrained local minimizer for sufficiently small $h$.
Theorem~\ref{thm:discrete-first-order-en} therefore recovers a unique
$p_h\in Q_h$ and the full system \eqref{eq:discrete-optimality-en}.

\begin{proposition}
\label{prop:coupled-error-estimates-en}
Under the assumptions of Theorem~\ref{thm:control-error-en}, the state and
adjoint energy errors satisfy
\begin{align}
  E_T
  &\leq
  ChM_{\vbar}|\Tbar|_2
  +
  C\kappa^{-2}\norm[0]{f}E_v,
  \label{eq:coupled-state-error-en}\\
  E_q
  &\leq
  ChM_{\vbar}|\qbar|_2
  +
  C\kappa^{-1}E_T
  +
  C\kappa^{-3}\norm[0]{f}E_v.
  \label{eq:coupled-adjoint-error-en}
\end{align}
and the pressure error satisfies
\begin{equation}
\label{eq:pressure-error-en}
\begin{aligned}
    E_p
  \leq{}&
  Ch|\pbar|_1
  +
  C\gamma\bigl(E_v+h|\vbar|_2\bigr)+
  C\norm[0]{f}
  \left(
  \kappa^{-2}E_T+\kappa^{-1}E_q
  \right).
\end{aligned}
\end{equation}
\end{proposition}

\begin{proof}
Using the auxiliary variables defined by
\eqref{eq:aux-state-en}--\eqref{eq:aux-adjoint-en}, we decompose the errors as
\begin{align*}
 \Tbar-T_h&=(\Tbar-\widehat T_h)+(\widehat T_h-T_h),\\
 \qbar-q_h&=(\qbar-\widehat q_h)+(\widehat q_h-q_h).
\end{align*}
For the state, the triangle inequality, the approximation estimate in
Lemma~\ref{lem:auxiliary-errors-en}, and the propagation estimate in
Lemma~\ref{lem:coupling-propagation-en}, with $\vh=\vbar_h$, imply
\begin{align*}
 E_T
 &\leq\norm[1]{\Tbar-\widehat T_h}
       +\norm[1]{\widehat T_h-T_h}\leq ChM_{\vbar}|\Tbar|_2
       +C\kappa^{-2}\norm[0]{f}E_v.
\end{align*}
For the adjoint, the same lemmas yield
\begin{align*}
 E_q
 &\leq\norm[1]{\qbar-\widehat q_h}
       +\norm[1]{\widehat q_h-q_h}\leq ChM_{\vbar}|\qbar|_2
       +C\kappa^{-1}E_T
       +C\kappa^{-3}\norm[0]{f}E_v.
\end{align*}
These inequalities prove
\eqref{eq:coupled-state-error-en}--\eqref{eq:coupled-adjoint-error-en}.

For the pressure, we decompose the error as
\[
 \pbar-p_h=(\pbar-P_h\pbar)+\eta_p,
 \qquad \eta_p:=P_h\pbar-p_h\in Q_h.
\]
The first term is controlled by \eqref{eq:pressure-projection-en}, while
the discrete inf--sup condition applies to $\eta_p$.
For every $\wh\in\Uh$, the inclusion $\nabla\cdot\wh\in Q_h$ and
the $L^2$ orthogonality of $P_h$ imply
\[
 b(\wh,\pbar-P_h\pbar)
 =\ip{\nabla\cdot\wh}{\pbar-P_h\pbar}=0.
\]
Consequently,
\[
 b(\wh,\pbar-p_h)
 =b(\wh,\pbar-P_h\pbar)+b(\wh,P_h\pbar-p_h)
 =b(\wh,\eta_p).
\]
Subtracting \eqref{eq:discrete-optimality-control-en} from
\eqref{eq:continuous-control-tested-en}, with the same test function
$\wh\in\Uh$, we obtain
\begin{align}
 b(\wh,\pbar-p_h)
 ={}&\gamma a_h(\vbar-\vbar_h,\wh)+\gamma\delta_h(\vbar,\wh)-\ip{\qbar\nabla\Tbar-q_h\nabla T_h}{\wh}.
 \label{eq:pressure-error-equation-en}
\end{align}
To estimate the nonlinear term, we write
\[
 \qbar\nabla\Tbar-q_h\nabla T_h
 =(\qbar-q_h)\nabla\Tbar+q_h\nabla(\Tbar-T_h).
\]
H\"older's inequality, the embedding
$H_0^1(\Omega)\hookrightarrow L^4(\Omega)$, and
Lemma~\ref{lem:discrete-trilinear-en} imply
\begin{align*}
 \left|\ip{\qbar\nabla\Tbar-q_h\nabla T_h}{\wh}\right|
 \leq{}&\norm[L^4(\Omega)]{\qbar-q_h}\norm[1]{\Tbar}
          \norm[L^4(\Omega)]{\wh}\\*
       &+\norm[L^4(\Omega)]{q_h}\norm[1]{\Tbar-T_h}
          \norm[L^4(\Omega)]{\wh}\\*
 \leq{}&C\bigl(E_q\norm[1]{\Tbar}
              +\norm[1]{q_h}E_T\bigr)\tnorm{\wh}.
\end{align*}
By \eqref{eq:state-stability-en} and
\eqref{eq:discrete-adjoint-stability-en}, we obtain
\[
 \left|\ip{\qbar\nabla\Tbar-q_h\nabla T_h}{\wh}\right|
 \leq C\norm[0]{f}\left(\kappa^{-1}E_q+\kappa^{-2}E_T\right)\tnorm{\wh}.
\]
For the other two terms in \eqref{eq:pressure-error-equation-en},
Cauchy--Schwarz for $a_h$ and \eqref{eq:consistency-bound-en} imply
\begin{align*}
 \gamma|a_h(\vbar-\vbar_h,\wh)|
 &\leq\gamma\tnorm{\vbar-\vbar_h}\tnorm{\wh}
  =\gamma E_v\tnorm{\wh},\\*
 \gamma|\delta_h(\vbar,\wh)|
 &\leq C\gamma h|\vbar|_2\tnorm{\wh}.
\end{align*}
Combining these bounds with
$b(\wh,\eta_p)=b(\wh,\pbar-p_h)$ and
\eqref{eq:pressure-error-equation-en}, we conclude that, for every
nonzero $\wh\in\Uh$,
\[
 \frac{|b(\wh,\eta_p)|}{\tnorm{\wh}}
 \leq C\gamma\bigl(E_v+h|\vbar|_2\bigr)
       +C\norm[0]{f}\bigl(\kappa^{-1}E_q+\kappa^{-2}E_T\bigr).
\]
The right-hand side is independent of $\wh$.  Taking the supremum and
using the discrete inf--sup constant $\beta_0>0$, we obtain
\begin{align*}
 \norm[0]{\eta_p}
 &\leq\frac1{\beta_0}
       \sup_{\bm0\ne\wh\in\Uh}
       \frac{|b(\wh,\eta_p)|}{\tnorm{\wh}}\\*
 &\leq\frac C{\beta_0}\left[
       \gamma\bigl(E_v+h|\vbar|_2\bigr)
       +\norm[0]{f}\bigl(\kappa^{-2}E_T+\kappa^{-1}E_q\bigr)\right].
\end{align*}
Finally, the triangle inequality and \eqref{eq:pressure-projection-en}
imply
\begin{align*}
 E_p=\norm[0]{\pbar-p_h}
 &\leq\norm[0]{\pbar-P_h\pbar}+\norm[0]{\eta_p}\\*
 &\leq Ch|\pbar|_1
       +\frac C{\beta_0}\left[
          \gamma\bigl(E_v+h|\vbar|_2\bigr)
          +\norm[0]{f}\bigl(\kappa^{-2}E_T+\kappa^{-1}E_q\bigr)
        \right].
\end{align*}
Absorbing the mesh-independent factor $\beta_0^{-1}$ into $C$ proves
\eqref{eq:pressure-error-en}.
\end{proof}

\begin{corollary}
\label{cor:full-error-en}
Under the assumptions of Theorem~\ref{thm:control-error-en},
the combined error satisfies
\begin{equation}
\label{eq:full-error-en}
  E_T+E_q+E_v+E_p
  \leq
  C_{\kappa,\gamma,\mu,\mathrm{data}}h.
\end{equation}
\end{corollary}

\begin{proof}
Theorem~\ref{thm:control-error-en} establishes $E_v\leq Ch$.
Successive substitution into \eqref{eq:coupled-state-error-en},
\eqref{eq:coupled-adjoint-error-en}, and \eqref{eq:pressure-error-en}
yields $E_T\leq Ch$, $E_q\leq Ch$, and $E_p\leq Ch$, respectively.
Adding these bounds proves \eqref{eq:full-error-en}; the constants depend
only on the fixed parameters and regularity quantities displayed above.
\end{proof}

\section{Numerical experiments}
\label{sec:numerics-en}

The computational domain is $\Omega=(0,1)^2$.  The state and adjoint are approximated using
continuous nodal $P_1$ basis functions associated with interior vertices.  The control is
approximated in $\Uh=P_1^c\oplus RT_0$, using the componentwise nodal
basis and the oriented flux basis $\{\bm\Phi_e\}_{e\in\Ehzero}$.
The pressure is piecewise constant
in $Q_h=P_0$.  The triangular mesh is generated using $n$ uniform subdivisions
along each coordinate direction.
We take $\alpha_K=1$ and use the diagonal form
$a_R=a_R^{\mathrm{div}}$.
For every computed solution we report
\[
  J_{\rm tr}:=\frac12\norm[0]{DT_h}^2,
  \qquad
  J_{\rm reg}:=\frac\gamma2a_h(\vh,\vh),
  \qquad
  j_h(\vh)=J_{\rm tr}+J_{\rm reg},
\]
and the relative reduction
$100\bigl(j_h(\bm0)-j_h(\vh)\bigr)/j_h(\bm0)$.  The factor $100$ converts the
relative decrease of the reduced functional into a percentage.  The divergence of control
\[
  E_{\rm div}:=\norm[0]{\nabla\cdot\vh}
\]
is evaluated elementwise from the full $P_1^c\oplus RT_0$ velocity.

\exampleheading
\label{subsec:manufactured-en}

The physical optimization problem has no known exact solution.  To measure
discretization error without conflating it with the selection of a local
solution, we prescribe the smooth fields
\begin{align*}
  T(x,y)&=\sin(\pi x)\sin(\pi y),
  &q(x,y)&=\sin(2\pi x)\sin(\pi y),\\
  \Psi(x,y)&=0.05\sin^2(\pi x)\sin^2(\pi y),
  &\bm v&=(\partial_y\Psi,-\partial_x\Psi),\\
  p(x,y)&=\sin(\pi x)\cos(\pi y),
  &\kappa&=1.
\end{align*}
The stream-function construction makes $\bm v$ exactly divergence-free and
ensures that it vanishes on $\partial\Omega$.  We define
\begin{align*}
 f&=-\kappa\Delta T+\bm v\cdot\nabla T,\\
 r&=-\kappa\Delta q-\bm v\cdot\nabla q-DT,\\
 \bm g_\gamma&=-\gamma\Delta\bm v+\nabla p-q\nabla T,
\end{align*}
and add $r$ and $\bm g_\gamma$ to the adjoint and control equations,
respectively.  Thus $(T,q,\bm v,p)$ solves the complete forced optimality
system exactly.

We use the errors defined in \eqref{eq:energy-error-notation-en}, with the
manufactured fields $(T,q,\bm v,p)$ as the exact solution.
Table~\ref{tab:convergence-energy-en} reports the
errors; parenthesized values are the observed orders computed from successive
values of $n$.

\begin{table}[!htbp]
\centering
\caption{Convergence rates for two regularization parameters.}
\label{tab:convergence-energy-en}
\small
\setlength{\tabcolsep}{1.25pt}
\begin{tabular}{@{}cccccc@{}}
\toprule
$n$ & $E_T$ & $E_q$ & $E_v$ & $E_p$ & $E_{\rm div}$ \\
\midrule
\multicolumn{6}{c}{$\gamma=1$} \\
\cmidrule(lr){1-6}
$4$   & $8.39\times10^{-1}$ & $1.89$ & $6.10\times10^{-1}$
      & $2.42\times10^{-1}$ & $4.24\times10^{-17}$ \\
$8$   & $4.32\times10^{-1}\;(0.96)$ & $1.00\;(0.91)$
      & $3.45\times10^{-1}\;(0.82)$ & $9.83\times10^{-2}\;(1.30)$
      & $8.13\times10^{-17}$ \\
$16$  & $2.18\times10^{-1}\;(0.99)$ & $5.09\times10^{-1}\;(0.98)$
      & $1.79\times10^{-1}\;(0.95)$ & $4.24\times10^{-2}\;(1.21)$
      & $1.40\times10^{-16}$ \\
$32$  & $1.09\times10^{-1}\;(1.00)$ & $2.55\times10^{-1}\;(0.99)$
      & $9.02\times10^{-2}\;(0.99)$ & $2.02\times10^{-2}\;(1.07)$
      & $2.62\times10^{-16}$ \\
$64$  & $5.45\times10^{-2}\;(1.00)$ & $1.28\times10^{-1}\;(1.00)$
      & $4.52\times10^{-2}\;(1.00)$ & $1.00\times10^{-2}\;(1.02)$
      & $6.50\times10^{-16}$ \\
$128$ & $2.73\times10^{-2}\;(1.00)$ & $6.39\times10^{-2}\;(1.00)$
      & $2.26\times10^{-2}\;(1.00)$ & $4.99\times10^{-3}\;(1.00)$
      & $1.04\times10^{-15}$ \\
\addlinespace[2pt]
\midrule
\multicolumn{6}{c}{$\gamma=10^{-2}$} \\
\cmidrule(lr){1-6}
$4$   & $8.39\times10^{-1}$ & $1.89$ & $2.37$
      & $1.62\times10^{-1}$ & $1.15\times10^{-16}$ \\
$8$   & $4.32\times10^{-1}\;(0.96)$ & $1.00\;(0.91)$
      & $7.20\times10^{-1}\;(1.72)$ & $7.18\times10^{-2}\;(1.17)$
      & $9.62\times10^{-17}$ \\
$16$  & $2.18\times10^{-1}\;(0.99)$ & $5.09\times10^{-1}\;(0.98)$
      & $2.41\times10^{-1}\;(1.58)$ & $3.37\times10^{-2}\;(1.09)$
      & $1.38\times10^{-16}$ \\
$32$  & $1.09\times10^{-1}\;(1.00)$ & $2.55\times10^{-1}\;(0.99)$
      & $9.88\times10^{-2}\;(1.28)$ & $1.65\times10^{-2}\;(1.03)$
      & $2.56\times10^{-16}$ \\
$64$  & $5.45\times10^{-2}\;(1.00)$ & $1.28\times10^{-1}\;(1.00)$
      & $4.64\times10^{-2}\;(1.09)$ & $8.20\times10^{-3}\;(1.01)$
      & $6.77\times10^{-16}$ \\
$128$ & $2.73\times10^{-2}\;(1.00)$ & $6.39\times10^{-2}\;(1.00)$
      & $2.28\times10^{-2}\;(1.03)$ & $4.09\times10^{-3}\;(1.00)$
      & $1.07\times10^{-15}$ \\
\bottomrule
\end{tabular}
\end{table}

\begin{figure}[!htbp]
  \centering
  \includegraphics[width=\textwidth]{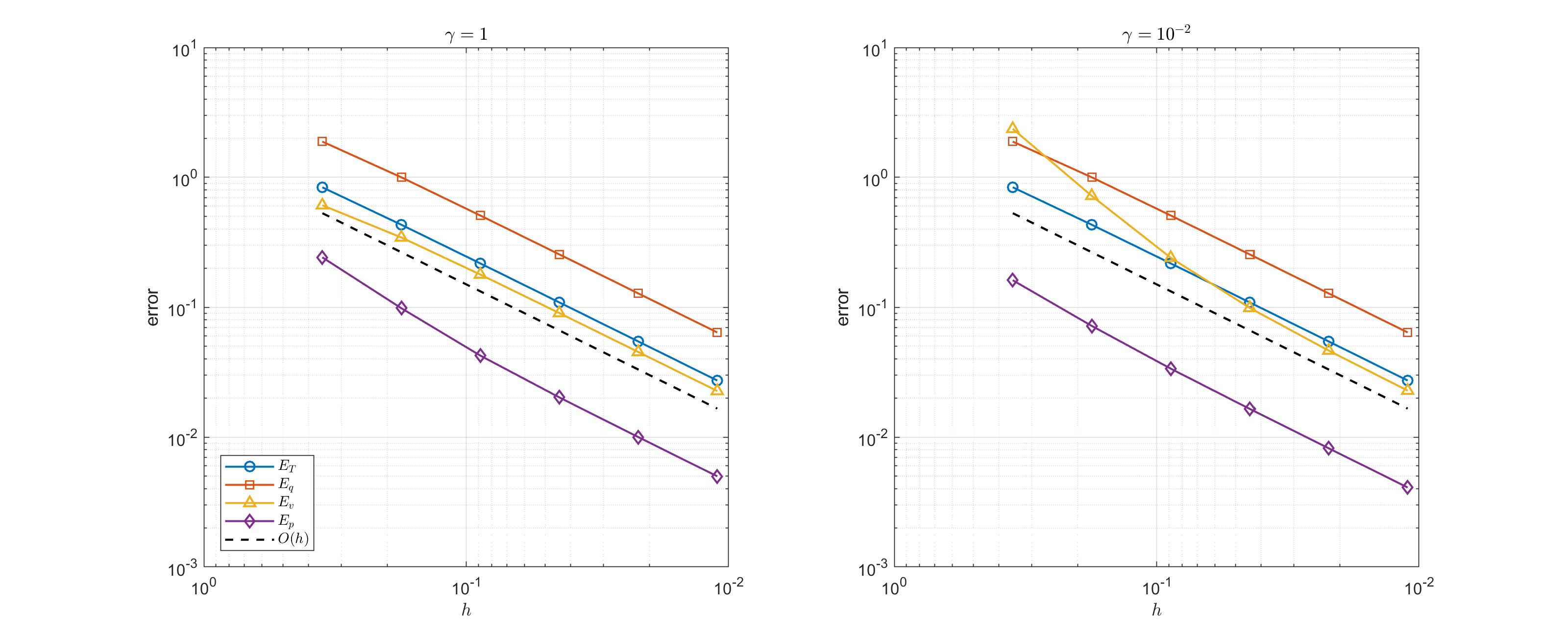}
  \caption{Log--log convergence histories.  Left: $\gamma=1$.  Right:
  $\gamma=10^{-2}$.  The reference lines have slope one.}
  \label{fig:convergence-rates-en}
\end{figure}

Table~\ref{tab:convergence-energy-en} and
Figure~\ref{fig:convergence-rates-en} show first-order convergence of all four
unknowns for both values of $\gamma$.  The control error for
$\gamma=10^{-2}$ is larger on coarse meshes and reaches its asymptotic regime
later, but its slope approaches one under refinement.  In both cases
$E_{\rm div}$ remains at the level of accumulated roundoff error, in agreement with the
pointwise mass-conservation property.

\exampleheading

We first consider the positive symmetric heat-source data from
Example~5.1 of \cite{HeHuMu2021-en}, namely
\[
  f(x,y)=2\pi^2\sin(\pi x)\sin(\pi y),
  \qquad \kappa=1.
\]
  The value of the uncontrolled tracking functional is $j_h(\bm0)=4.287\times10^{-2}$.
Table~\ref{tab:reference-case1-en} shows that decreasing $\gamma$ produces a
stronger recirculating velocity, reduces the temperature maximum, and lowers
the objective.

\begin{table}[htbp]
\centering
\caption{Positive-source test problem on the $n=100$ mesh. The velocity maximum is evaluated at the triangular
quadrature points.}
\label{tab:reference-case1-en}
\footnotesize
\setlength{\tabcolsep}{1.7pt}
\begin{tabular}{@{}cccccccc@{}}
\toprule
$\gamma$ & $J_{\rm tr}$ & $J_{\rm reg}$ & $j_h$ & reduction (\%)
& $\max T_h$ & $\max|\vh|$ & $E_{\rm div}$ \\
\midrule
$3.60\times10^{-6}$ & $4.193\times10^{-2}$ & $7.023\times10^{-4}$
& $4.263\times10^{-2}$ & $0.54$ & $0.9926$ & $2.3866$ & $3.86\times10^{-14}$ \\
$8.50\times10^{-7}$ & $2.355\times10^{-2}$ & $1.001\times10^{-2}$
& $3.355\times10^{-2}$ & $21.73$ & $0.7763$ & $18.2080$ & $3.00\times10^{-13}$ \\
$3.90\times10^{-7}$ & $1.648\times10^{-2}$ & $9.375\times10^{-3}$
& $2.586\times10^{-2}$ & $39.67$ & $0.6654$ & $25.6653$ & $4.16\times10^{-13}$ \\
\bottomrule
\end{tabular}
\end{table}
\begin{figure}[!htbp]
  \centering
  \includegraphics[width=0.8\textwidth]
  {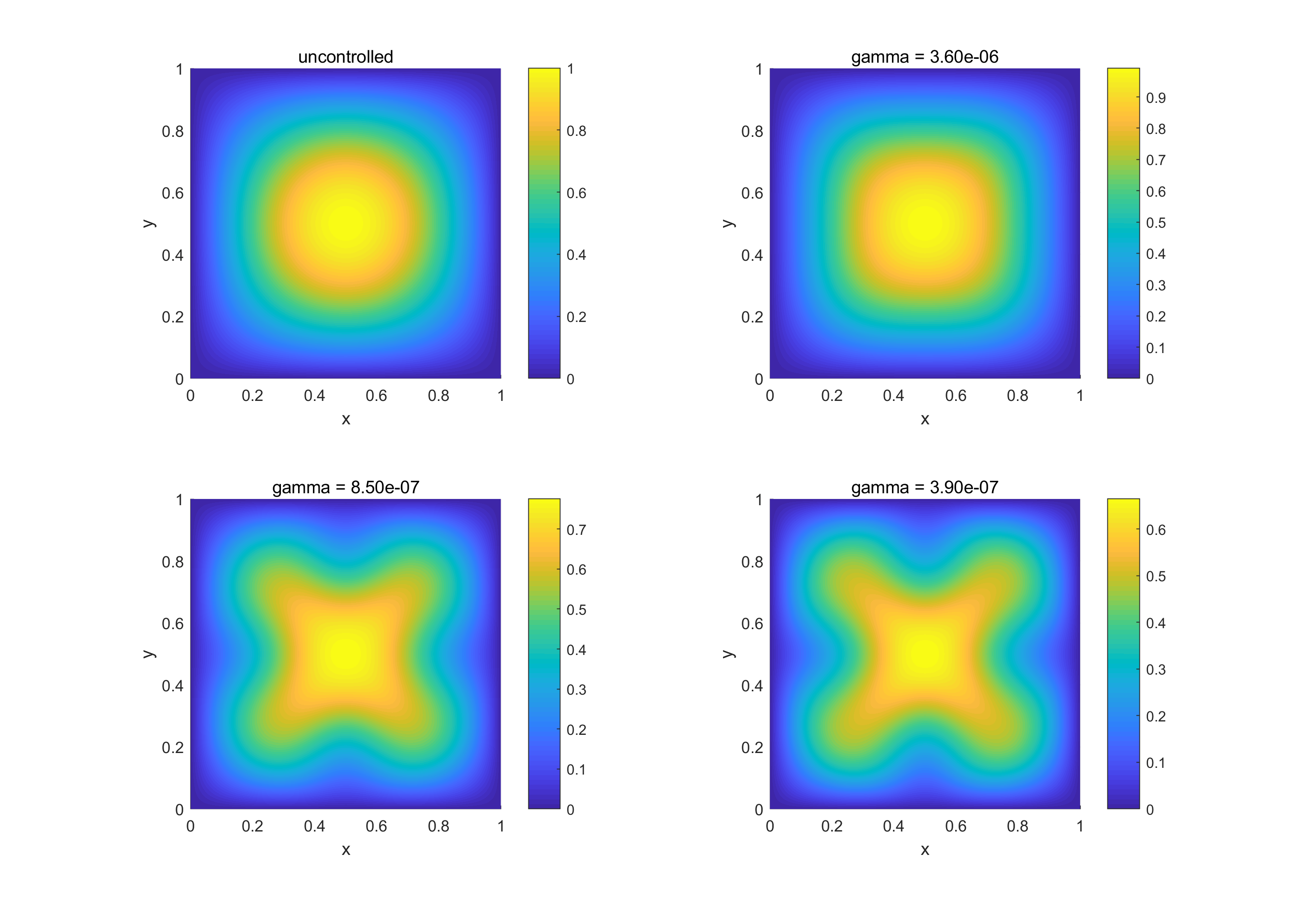}
  \caption{Temperature for the positive-source test problem: uncontrolled field
  and three decreasing values of $\gamma$ on the $n=100$ mesh.}
  \label{fig:reference-case1-temperature-en}
\end{figure}
Figure~\ref{fig:reference-case1-temperature-en} shows the corresponding
change in the spatial temperature distribution.  With no control, the temperature has a
single smooth maximum at the centre and nearly concentric level sets.  The
weakly controlled field at $\gamma=3.6\times10^{-6}$ remains close to this
diffusion-dominated profile.  As $\gamma$ decreases, the stronger velocity field
deforms the central hot region into a fourfold-symmetric pattern and spreads
the heat along four preferred transport directions.  This redistribution is
consistent with the simultaneous decrease of $J_{\rm tr}$ and $\max T_h$ in
Table~\ref{tab:reference-case1-en}.

\begin{figure}[!htbp]
  \centering
  \includegraphics[width=0.8\textwidth]
  {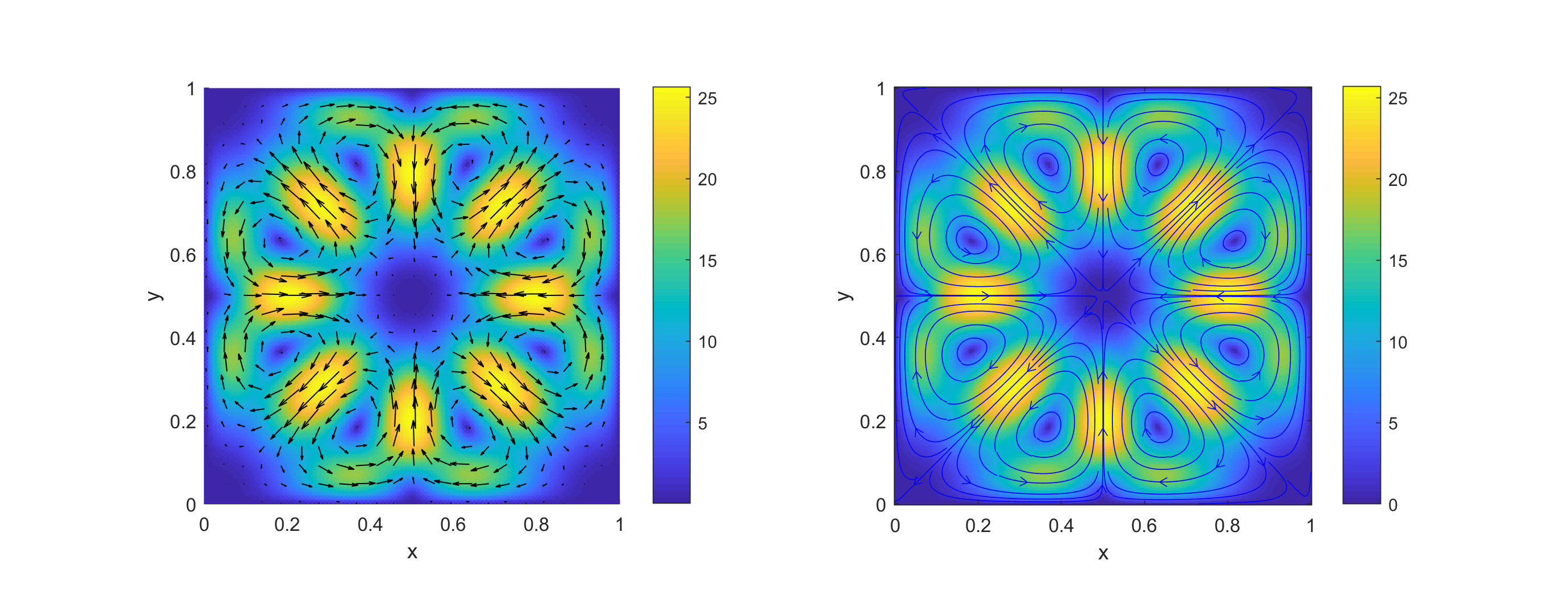}
  \caption{Computed velocity for the positive-source test problem at
  $\gamma=3.9\times10^{-7}$ on the $n=100$ mesh.  Left: the raw
  $P_1^c\oplus RT_0$ velocity at element centroids.  Right: streamlines of its
  continuous $L^2$ projection used only for visualization. }
  \label{fig:reference-case1-velocity-en}
\end{figure}

The flow mechanism responsible for this redistribution is visible in
Figure~\ref{fig:reference-case1-velocity-en}.  Eight recirculation cells are
arranged around the centre of the square, and the largest speeds occur in the
interior bands separating neighbouring cells.  Their alternating circulation
transports heat away from the central peak while producing no net boundary
flux.

To assess mesh dependence, we repeat the computation at the last
parameter value on meshes with $n=32,64,100$.  Table~\ref{tab:reference-case1-mesh-en} shows only small
changes in the temperature maximum, velocity maximum, and objective between
$n=64$ and $n=100$, indicating that these scalar outputs are nearly
mesh-independent at the two finest resolutions.

\begin{table}[htbp]
\centering
\caption{Mesh-refinement study for the positive-source test problem at
$\gamma=3.9\times10^{-7}$.}
\label{tab:reference-case1-mesh-en}
\small
\begin{tabular}{@{}ccccccc@{}}
\toprule
$n$ & triangles & $\max T_h$ & $\max|\vh|$ & $j_h$ & reduction (\%)
& $E_{\rm div}$ \\
\midrule
$32$  & $4\,096$  & $0.666711$ & $25.8219$ & $2.587\times10^{-2}$
& $39.588$ & $6.26\times10^{-14}$ \\
$64$  & $16\,384$ & $0.665603$ & $25.7006$ & $2.586\times10^{-2}$
& $39.660$ & $1.22\times10^{-13}$ \\
$100$ & $40\,000$ & $0.665385$ & $25.6653$ & $2.586\times10^{-2}$
& $39.674$ & $4.16\times10^{-13}$ \\
\bottomrule
\end{tabular}
\end{table}

\FloatBarrier

\exampleheading

To test sign-changing data and nonradial transport, we next consider the
source--sink data of Example~5.5 in \cite{HeHuMu2021-en}:
\[
  f(x,y)=75\exp\!\left(-\frac{(9x-2)^2+(9y-2)^2}{4}\right)
       -75\exp\!\left(-\frac{(9x-4)^2+(9y-7)^2}{4}\right)
\]
with $\kappa=1$.  We first use the $n=32$ mesh
($4\,096$ triangles), for which
$j_h(\bm0)=1.289\times10^{-1}$.
Table~\ref{tab:reference-case5-en} reports the results for three regularization
parameters.  Decreasing $\gamma$ permits a much stronger velocity: the maximum
speed rises from $2.74$ to $22.69$.
The very small values of
$E_{\rm div}$ show that this enhanced transport is generated without losing
discrete mass conservation.

\begin{table}[htbp]
\centering
\caption{Source--sink test problem on the $n=32$ mesh.}
\label{tab:reference-case5-en}
\small
\resizebox{\textwidth}{!}{%
\begin{tabular}{@{}ccccccccc@{}}
\toprule
$\gamma$ & $J_{\rm tr}$ & $J_{\rm reg}$ & $j_h$ & reduction (\%)
& $\min T_h$ & $\max T_h$ & $\max|\vh|$ & $E_{\rm div}$ \\
\midrule
$5.00\times10^{-5}$ & $1.257\times10^{-1}$ & $2.805\times10^{-3}$
& $1.285\times10^{-1}$ & $0.321$ & $-1.3034$ & $1.0487$ & $2.7439$ & $6.59\times10^{-15}$ \\
$1.00\times10^{-5}$ & $6.208\times10^{-2}$ & $2.937\times10^{-2}$
& $9.146\times10^{-2}$ & $29.062$ & $-0.8178$ & $0.7601$ & $19.3924$ & $4.91\times10^{-14}$ \\
$6.90\times10^{-6}$ & $5.278\times10^{-2}$ & $2.800\times10^{-2}$
& $8.078\times10^{-2}$ & $37.342$ & $-0.7446$ & $0.6960$ & $22.6904$ & $6.09\times10^{-14}$ \\
\bottomrule
\end{tabular}}
\end{table}

\begin{figure}[!htbp]
  \centering
  \includegraphics[width=0.8\textwidth]{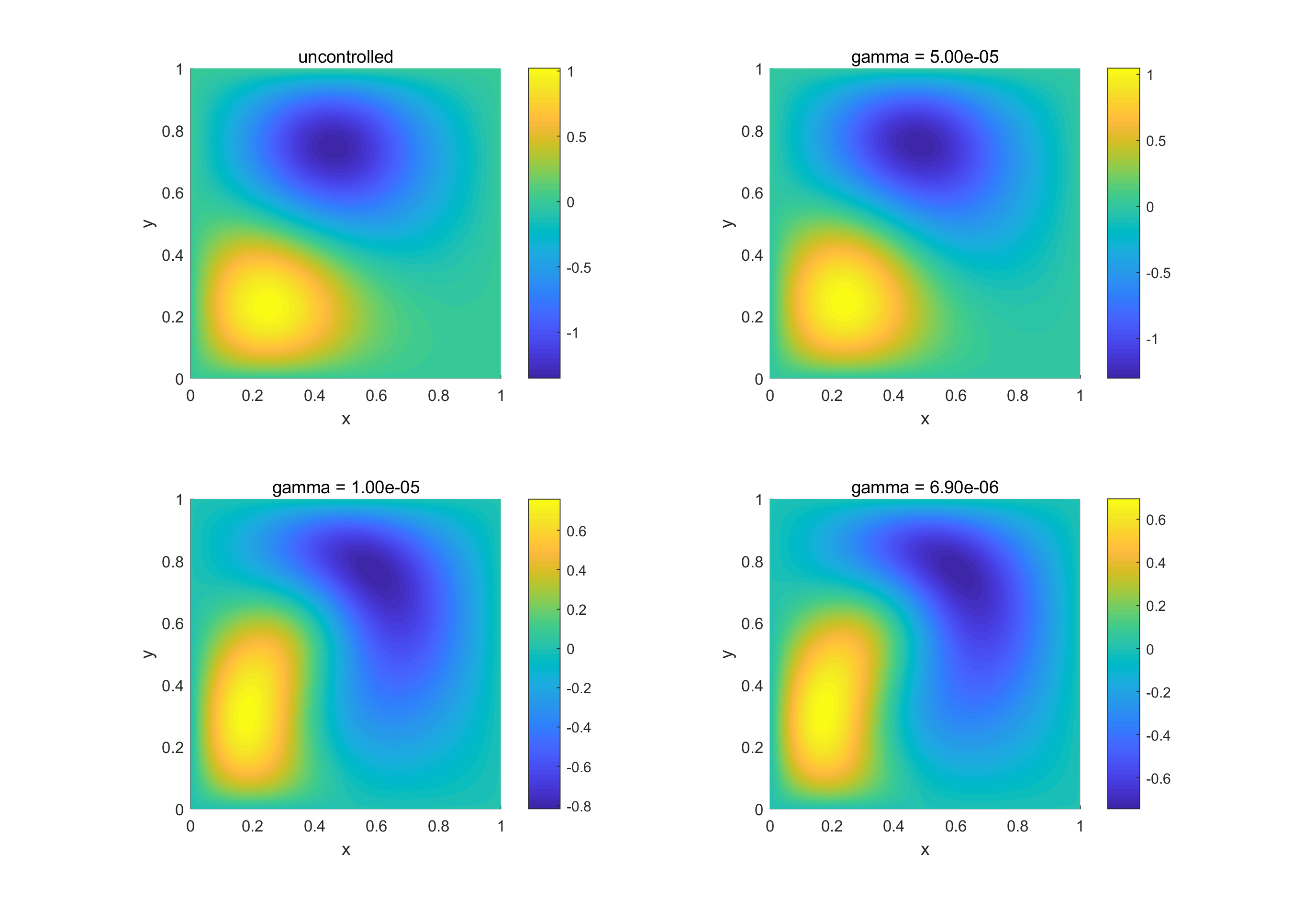}
  \caption{Temperature for the source--sink test problem: uncontrolled field
  and three decreasing values of $\gamma$ on the $n=32$ mesh.}
  \label{fig:reference-case5-temperature-en}
\end{figure}

Figure~\ref{fig:reference-case5-temperature-en} shows that the control changes
both the amplitude and the transport path of the temperature.  For
$\gamma=5\times10^{-5}$ the velocity penalty dominates, and the temperature
remains close to the uncontrolled field with a hot region in the lower-left
part of the square and a cold region above it.  At $\gamma=10^{-5}$ the two
regions are stretched and separated by convection.  For the selected discrete
solution at $\gamma=6.9\times10^{-6}$, the warm region is elongated upward
along the left part of the domain, whereas the cold region spreads through
the upper-right part.

\begin{figure}[!htbp]
  \centering
  \includegraphics[width=\textwidth]{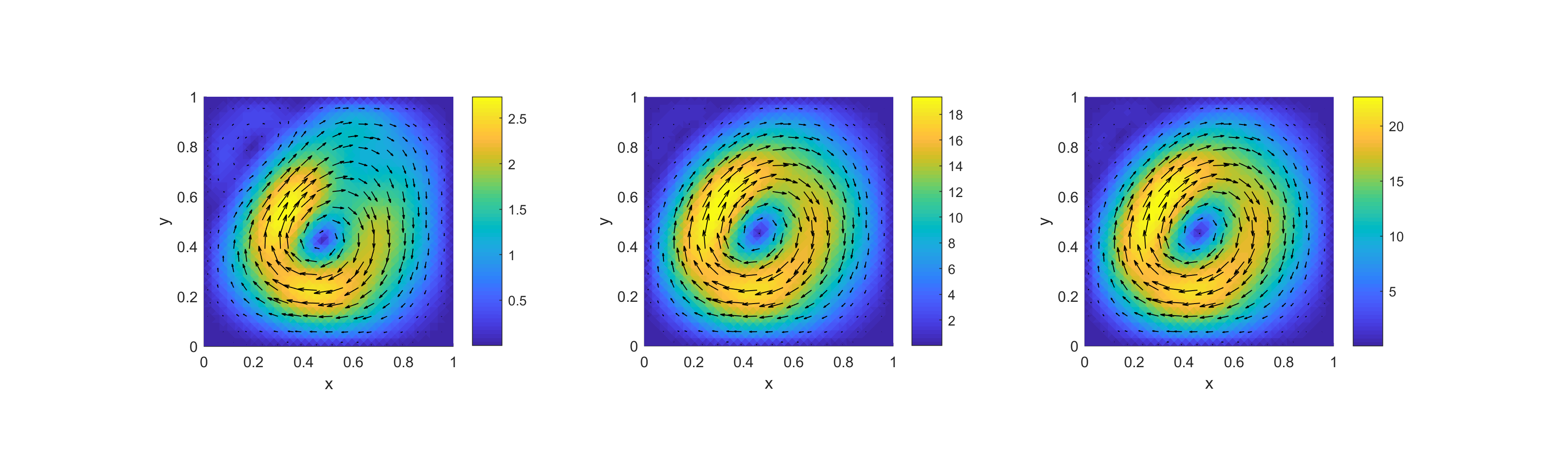}
  \caption{Actual $P_1^c\oplus RT_0$ controls for the source--sink problem,
  evaluated at element centroids on the $n=32$ mesh.  From left to right:
  $\gamma=5.0\times10^{-5}$, $10^{-5}$, and $6.9\times10^{-6}$.}
  \label{fig:reference-case5-velocity-en}
\end{figure}
Figure~\ref{fig:reference-case5-velocity-en} shows how weaker regularization
changes the strength of the optimal control.  At
$\gamma=5.0\times10^{-5}$ the velocity is weak and forms one broad
recirculation cell.  The same cell becomes much stronger at $\gamma=10^{-5}$.
At $\gamma=6.9\times10^{-6}$, the single clockwise circulation cell
strengthens further, with the highest speeds surrounding the interior
low-speed region.
This flow pattern is consistent with the elongated temperature
distribution in Figure~\ref{fig:reference-case5-temperature-en}.
The velocity magnitude decays toward the no-slip boundary in all three cases.

\exampleheading
\label{subsec:pressure-robustness-test-en}

The divergence diagnostic $E_{\rm div}$ verifies pointwise mass conservation,
but it does not by itself test pressure robustness.  We use the smooth forced
optimality system from Example~\ref{subsec:manufactured-en}.  To isolate
the effect of an irrotational force, we perturb the pressure and control forcing
as follows:
\[
  \psi(x,y)=\sin(\pi x)\sin(\pi y),
  \qquad
  p_\Lambda=p+\Lambda\bigl(\psi-\langle\psi\rangle\bigr),
  \qquad
  \bm g_\Lambda=\bm g+\Lambda\nabla\psi,
\]
where $\bm g=\bm g_\gamma$ is the original manufactured control forcing and
\[
  \Lambda\in\{1,10^2,10^4,10^6,10^8\}.
\]
The two terms are changed together so that
\[
  -\gamma\Delta\bm v+\nabla p_\Lambda
  =q\nabla T+\bm g_\Lambda.
\]
Thus the exact state $T$, adjoint $q$, and velocity $\bm v$ are independent of
$\Lambda$; only the pressure amplitude changes.  Any variation in the
computed velocity measures pressure sensitivity when the nonlinear solve
selects the corresponding solution and quadrature and algebraic errors are controlled.
Exact divergence freedom and zero normal trace imply
$\ip{\nabla\psi}{\zh}=0$ for all $\zh\in\Zh$ and $\psi\in H^1(\Omega)$.
With exact integration, adding $\Lambda\nabla\psi$ to the control forcing
leaves $T_h$, $q_h$, and $\vh$ unchanged and shifts the pressure to
$p_h+\Lambda P_h\psi$, where $P_h$ is the $L^2$-orthogonal projection
onto the mean-zero space $Q_h$.

For comparison with a weakly divergence-free method, we solve the same complete forced
optimality system with the classical MINI velocity and continuous $P_1$
pressure pair. The velocity space consists of continuous $P_1$ functions enriched
by one cubic bubble per component on each triangle.  The state and adjoint
remain continuous $P_1$ functions, and a skew-symmetric convection form is
used for both methods.  Thus the comparison changes only the
velocity--pressure discretization and does not introduce a different energy
stabilization.
For the MINI method, the control equation is derived consistently
from the discrete optimization problem with skew-symmetric convection,
yielding $\frac12(q_h\nabla T_h-T_h\nabla q_h,\bm w_h)$
as the convection contribution.
We fix $n=16$, $\kappa=1$, and $\gamma=10^3$, and vary
$\Lambda$ from $1$ to $10^8$, with $\Lambda=0$ as an additional baseline.
Table~\ref{tab:pressure-robustness-en} reports the natural velocity energy
error $E_v$, the velocity error
$E_{v,0}:=\norm[0]{\bm v-\bm v_h}$, and $E_{\rm div}$.  For the MINI pair,
$E_v$ is the $H^1$ seminorm of the velocity error.

\begin{table}[!htbp]
\centering
\caption{Gradient-force test on the fixed $n=16$ mesh with
$\kappa=1$ and $\gamma=10^3$.}
\label{tab:pressure-robustness-en}
\small
\begin{tabular}{@{}llrrr@{}}
\toprule
$\Lambda$ & method & $E_v$ & $E_{v,0}$
& $E_{\rm div}$ \\
\midrule
$0$ & $P_1^c\oplus RT_0-P_0$
& $1.786\times10^{-1}$ & $2.597\times10^{-3}$
& $1.23\times10^{-16}$ \\
& MINI--$P_1$
& $1.057\times10^{-1}$ & $2.571\times10^{-3}$
& $7.34\times10^{-2}$ \\
\addlinespace
$1$ & $P_1^c\oplus RT_0-P_0$
& $1.786\times10^{-1}$ & $2.597\times10^{-3}$
& $1.24\times10^{-16}$ \\
& MINI--$P_1$
& $1.057\times10^{-1}$ & $2.571\times10^{-3}$
& $7.34\times10^{-2}$ \\
\addlinespace
$10^2$ & $P_1^c\oplus RT_0-P_0$
& $1.786\times10^{-1}$ & $2.597\times10^{-3}$
& $1.27\times10^{-16}$ \\
& MINI--$P_1$
& $1.057\times10^{-1}$ & $2.571\times10^{-3}$
& $7.34\times10^{-2}$ \\
\addlinespace
$10^4$ & $P_1^c\oplus RT_0-P_0$
& $1.786\times10^{-1}$ & $2.597\times10^{-3}$
& $1.17\times10^{-16}$ \\
& MINI--$P_1$
& $1.063\times10^{-1}$ & $2.573\times10^{-3}$
& $7.41\times10^{-2}$ \\
\addlinespace
$10^6$ & $P_1^c\oplus RT_0-P_0$
& $1.786\times10^{-1}$ & $2.597\times10^{-3}$
& $1.37\times10^{-16}$ \\
& MINI--$P_1$
& $1.137$ & $9.990\times10^{-3}$
& $9.81\times10^{-1}$ \\
\addlinespace
$10^8$ & $P_1^c\oplus RT_0-P_0$
& $1.786\times10^{-1}$ & $2.597\times10^{-3}$
& $1.20\times10^{-16}$ \\
& MINI--$P_1$
& $1.132\times10^{2}$ & $9.654\times10^{-1}$
& $9.783\times10^{1}$ \\
\bottomrule
\end{tabular}
\end{table}

\begin{figure}[!htbp]
   \centering
  \includegraphics[width=\textwidth]
  {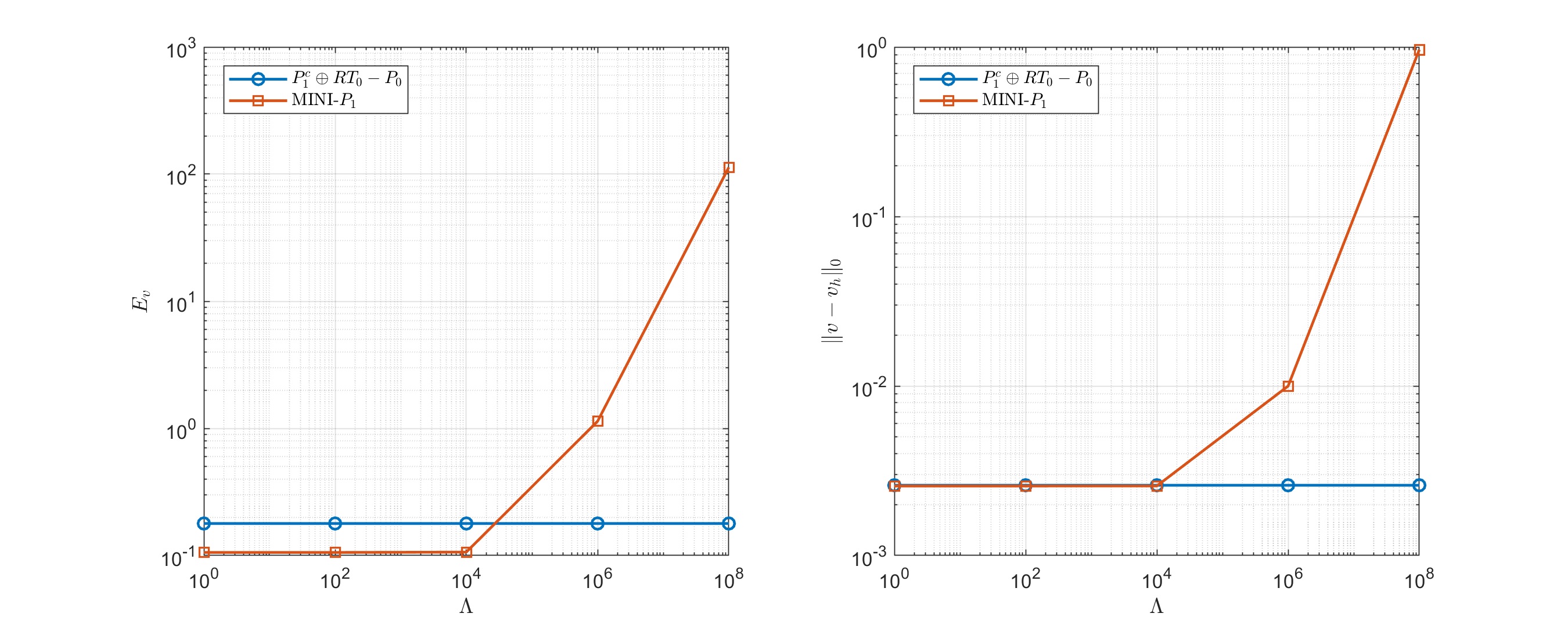}
  \caption{Velocity sensitivity to the gradient-force amplitude on a fixed
  mesh.  Left: energy error.  Right: $L^2$ error.}
  \label{fig:pressure-robustness-en}
\end{figure}

Figure~\ref{fig:pressure-robustness-en} shows that both velocity errors of
the proposed method remain independent of $\Lambda$ to the reported
precision, in agreement with the gradient-force invariance described above.  In contrast, the MINI
energy error increases from $1.057\times10^{-1}$ to $1.132\times10^2$, a
factor of approximately $1.07\times10^3$, and its $L^2$ error increases by a
factor of approximately $3.75\times10^2$. Relative to the $\Lambda=0$
baseline, the proposed velocity changes by less than $2.2\times10^{-12}$
in the discrete energy norm over the tested amplitudes. Finally, the value of
$E_{\rm div}$ for the proposed pair confirms mass conservation.

\FloatBarrier
\section*{\texorpdfstring{Conclusion}{Conclusion}}

The proposed discretization combines continuous $P_1$ elements for the
state and adjoint with a $P_1^c\oplus RT_0-P_0$ discretization of
the Stokes system satisfied by the optimal velocity. The velocity is
pointwise divergence-free and its approximation error is pressure robust.
Because the reduced problem is nonconvex, the analysis is performed for a
family of discrete local minimizers associated with an isolated continuous
strict local minimizer. The divergence-preserving interpolant, auxiliary
discrete variables, nonconforming consistency estimate, and reduced-gradient
residual relation together yield first-order a priori estimates for the
complete optimality system. The numerical results support the theoretical analysis.
In the future work we will consider the convection dominated case, where stablization techniques and adaptive refinement are needed to improve the numerical solutions.

\section*{Data Availability}
No data was used for the research described in the article.

\section*{Declarations}
The authors declare that they have no conflict of interest.

\bibliographystyle{siamplain}
\bibliography{divfree_optimal_control_BIB}

\end{document}